\documentclass[final,1p]{elsarticle}

\usepackage{graphicx}
\usepackage{subfigure}
\usepackage{epstopdf} 
\usepackage{amssymb,latexsym,amsmath,amsthm}
\biboptions{numbers,sort&compress}
\usepackage{hyperref}
\usepackage{bm,color}
\usepackage{multirow}
\usepackage{tabularx,booktabs}
\usepackage{stmaryrd}
\usepackage{lineno}
\usepackage{booktabs}       
\usepackage{float}
\usepackage{siunitx}        
\usepackage{array}          
\usepackage{caption}        
\usepackage{enumitem}
\usepackage{xcolor}
\numberwithin{equation}{section}
\allowdisplaybreaks

\newcommand{\f}{\frac}
\newcommand{\p}{\partial}

\newcommand{\mrm}{\mathrm}

\newcommand{\mal}{\mathcal}
\newcommand{\mA}{\mathcal{A}}
\newcommand{\mO}{\mathcal{O}}
\newcommand{\mT}{\mathcal{T}}
\newcommand{\mS}{\mathcal{S}}

\newcommand{\mE}{\mathcal{E}}
\newcommand{\mM}{\mathcal{M}}
\newcommand{\mK}{\mathcal{K}}
\newcommand{\mbb}{\mathbb}

\newcommand{\mfm}{\mathfrak{m}}

\DeclareMathOperator{\sech}{sech} 
\newtheorem{theorem}{Theorem}[section]
\newtheorem{lemma}[theorem]{Lemma}

\theoremstyle{definition}

\theoremstyle{remark}
\newtheorem{remark}{Remark}
\journal{XXX}
\begin{document}
   \begin{frontmatter}
      \title{Unconditional optimal-order error estimates of linear relaxation compact difference scheme for the coupled nonlinear Schr\"{o}dinger system}
		
	 \author[OUC]{Ying Gao}\ead{gy9168@stu.ouc.edu.cn}
	 \author[OUC,LMS]{Hongfei Fu\corref{Fu}} \ead{fhf@ouc.edu.cn}
     \author[OUC]{Xiaoying Wang} \ead{wxy7121@stu.ouc.edu.cn}

	 \address[OUC]{School of Mathematical Sciences, Ocean University of China, Qingdao, Shandong 266100, China}
	 \address[LMS]{Laboratory of Marine Mathematics, Ocean University of China, Qingdao, Shandong 266100, China}
	 \cortext[Fu]{Corresponding author.}
		
\begin{abstract}
    This paper presents a linear, decoupled, mass- and energy-conserving numerical scheme for the multi-dimensional coupled nonlinear Schr\"{o}dinger (CNLS) system. The scheme combines the fourth-order compact difference approximation in space with the relaxation technique in a time-staggered mesh framework, solving the primal unknowns and introduced auxiliary relaxation variables sequentially with high efficiency and high-order accuracy. We show the unique solvability and discrete conservation laws of the developed scheme. In particular, for the first time, leveraging an auxiliary error equation system combined with the cut-off technique, optimal-order error estimates in the discrete $H^1$-norm for the primal variables at the time nodes, and in the discrete $L^2$-norm for the auxiliary relaxation variables at the intermediate time nodes, are rigorously proved without any coupling mesh conditions, which contribute to the primary theoretical contribution of this paper for multi-dimensional CNLS system. Numerical experiments demonstrate convincingly the strong performance of the proposed scheme in long-term simulations, maintaining both physical invariants and high-order accuracy.
\end{abstract}
		
\end{frontmatter}

  \section{Introduction}
 In this paper, we develop and rigorously analyze a high-order fully discrete linearized numerical scheme for the coupled nonlinear Schr\"{o}dinger (CNLS) system
  \begin{equation}\label{model:e1a}
    \left\{
		\begin{aligned}
			& \mrm{i} u_{t}+\kappa \Delta u+ \left(|u|^{2}+\beta|v|^{2}\right)u=0, && \qquad \text{in} ~\Omega\times I,\\
			& \mrm{i} v_{t}+\kappa \Delta v+ \left(|v|^{2}+\beta|u|^{2}\right)v=0, && \qquad \text{in} ~\Omega\times I,
		\end{aligned}
        \right.
 \end{equation}
subject to the periodic boundary conditions and the following initial conditions
	\begin{equation}\label{model:e1b}
		u=u^{0},~~ v=v^{0}, \qquad \text{in} ~\Omega \times\{0\}, 
	\end{equation}
where $\Omega \subset \mbb{R}^d$ is an interval ($d=1$), or a rectangle ($d = 2$) or a cuboid ($d = 3$), the time interval $I:=(0,T]$ with $T< \infty$, $\kappa$ is the dispersion coefficient in the optical fiber, and $\beta$ is defined as a coupling parameter for birefringent optical fiber. The CNLS system couples multiple wavefields through two different nonlinear effects: self phase modulation (in field) and cross phase modulation (between fields), and it plays a critical role in such as quantum mechanics, nonlinear optics, and multicomponent Bose-Einstein condensates \cite{aydin2009,kanna2001,YWt2019}. From a numerical perspective, the CNLS system not only inherits the computational challenges of the nonlinear Schr\"{o}dinger equation (NLS), but also introduces new difficulties arising from the coupling of variables. There are a large number of literature on numerical methods and simulations of both the CNLS and NLS equations, to just name a few, finite difference methods \cite{AG1993,CQS1999,XSS2009,SWW2017,LHL2021,GZ2011},  finite element methods \cite{WJL2014,MS2008,GYL2016, MYH2022},  spectral methods \cite{BWZ2003,HMR1988} and discontinuous
Galerkin method \cite{KO1999,XY2005}. 

It is well known that the solutions of \eqref{model:e1a}--\eqref{model:e1b} possess the following \textit{mass conservation law}
	\begin{equation}\label{model:conserveMass}
		\mal{M}_1(t)=\int_{\Omega}|u|^2 \mrm{d}x=\mal{M}_1(0), \quad  \mal{M}_2(t)=\int_{\Omega}|v|^2 \mrm{d}x=\mal{M}_2(0),
	\end{equation}
	and \textit{energy conservation law}
	\begin{equation}\label{model:conserveEnergy}
		\mal{E}(t)  =\int_{\Omega}\Big[  \kappa  \left(|\nabla u|^{2}+ |\nabla v |^{2}\right)
		- \f{1}{2} \left(|u|^{4}+|v|^{4} \right)
		-  \beta  |u|^{2}|v|^{2}\Big]\mrm{d}x  
		=\mal{E}(0). 
	\end{equation}
Thus, the development of numerical methods that can preserve these physics-conserving properties is of great importance for long-term stable numerical simulations, and this has been one of the core research focuses for the numerical modeling of the CNLS equations. To date, extensive research has focused on developing mass- and/or energy-conserving numerical schemes for the NLS type equations, including the scalar auxiliary variable method \cite{BGM2024, FXB2021, LDF2023optimal, MS2023, GSM2023, GSM2025},  the Lagrange multiplier method \cite{GSM2023, GSM2025}, the integrating factor method \cite{ML2024,GXL2023,JCL2022}, the invariant energy quadratization method \cite{YY2021},  the quadratic auxiliary variable method \cite{FYY2023,HYF2023}, the multi-symplectic method \cite{HJL2006,HJL2009}, the operator splitting method \cite{LWY2015,ZSY2019}, the Hamiltonian boundary value method \cite{LB2018}, the average vector field method \cite{BJJ2023,AC2013}, and the fully-implicit discretization for nonlinear potentials \cite{GYZ2017,BWZ2012}.%

As another structure-preserving algorithm, the relaxation method is proposed by Besse for the NLS equation in \cite{BC1998,BC2004}, by introducing a new auxiliary relaxation variable, and then approximating the primal and auxiliary variables on staggered temporal grids via the Crank-Nicolson temporal discretization. Such method can achieve linearization and fully decoupling while rigorously preserving both mass and energy at the discrete level. However, the author only proves the convergence for the time-discrete scheme over a finite-time interval $[0,T]$ in \cite{BC2004}, while leaving the second-order temporal convergence analysis with respect to the temporal step size in \cite{BC2021} in the case of a cubic NLS equation.  Very recently, building upon this foundation, for the first time, Zouraris \cite{ZG2023} addresses the unconditional optimal-order error estimates for the fully discrete method that combines the relaxation technique with the spatial central finite difference approximation for the one-dimensional NLS equation. 
Further, Zhou et al. \cite{ZHOU2025} propose a fully discrete relaxation fourth-order compact finite difference method for the one-dimensional NLS equation and perform the corresponding optimal-order error estimates. However, the error analysis presented in \cite{ZG2023,ZHOU2025} critically depends on the discrete Sobolev inequality 
\begin{equation}\label{s1:e1}
  \|v\|_\infty\le C_{\Omega}\|v\|_1,~~\Omega \subset \mbb{R}, 
\end{equation}  
which directly provides a uniform bound for the numerical solution. Therefore, such unconditional error analysis for the fully discrete scheme is only limited to the one space dimension, as the discrete Sobolev inequality \eqref{s1:e1} does not hold for the multi-dimensional settings ($d=2, 3$). Consequently, as highlighted in \cite{ZG2023,ZHOU2025}, the techniques used to derive the unconditional optimal-order error estimates for the relaxation method in one-dimensional space 
cannot be directly generalized to the multi-dimensional space without imposing any coupling mesh condition.  
Therefore, there is currently a big gap on unconditional error analysis of the fully discrete relaxation numerical method for the multi-dimensional NLS type equations.

This paper aims to develop a linear, relaxation-based, compact difference (LRCD) scheme for the CNLS system, delivering high efficiency, high accuracy, and strict preservation of mass and energy at the discrete level. More importantly, for the first time, we establish a rigorous \textit{unconditional} and \textit{optimal} error analysis for the proposed scheme in multiple spatial dimensions up to $d= 3$. The main contributions of this paper are summarized as follows:
\begin{itemize}
    \item We rigorously prove that the LRCD scheme is uniquely solvable, by showing that the kernel of the corresponding linear operator is $\{0\}$, see Theorem \ref{thm:unique}, which guarantees the well-posedness of the numerical scheme.
    
    \item We demonstrate that the developed LRCD scheme strictly preserves mass and energy at the discrete level, see Theorem \ref{thm:conservation}, making it particularly suitable for long-term stable simulations.
    
    \item  We establish unconditional and optimal-order error estimates for the LRCD scheme in multiple space dimensions $d \le 3$, see Theorem \ref{thm:Convergence}). To the best of our knowledge, this is the first rigorous unconditional and optimal analysis for the relaxation-based difference scheme in multiple spatial dimensions. It is also applicable to relaxation methods combined with other spatial discretization techniques (e.g., finite element methods).
\end{itemize}
The unconditional and optimal-order error estimates of the LRCD shceme turns out to be very challenging, primarily due to the difficulty in proving the discrete $L^\infty$-norm bounds for numerical solutions without any coupling mesh condition \cite{ZG2023,ZHOU2025}.  To address this issue, in Section \ref{sec:AuEE} we introduce an auxiliary error equation system with four cut-off coefficients, in which a to-be-determined positive constant $\mK^*$ is involved in the  cut-off coefficients for the relaxation variables, see \eqref{cut_off:e1}. 
In Section \ref{sec:err_est}, through the time difference technique and the boundedness of the cut-off coefficients (see \eqref{ieq:bound_UmVm}--\eqref{ieq:bound_phi}), we first establish the unconditional optimal-order error estimates for the auxiliary error equation system in the discrete $L^2$-norm (see Theorem \ref{thm:M_err5}). Afterwards, to derive the discrete $L^\infty$-norm bounds of the numerical solutions to the auxiliary error equation system, we conduct a detailed discussion regarding the relationship between the temporal step size $\tau$ and spatial grid size $h$ (see Theorem \ref{thm:M_err6}):
(i) For $\tau>h^{d/2}$, our analysis proceeds in three stages. First, we estimate the discrete $H^2$-seminorm of the auxiliary errors for the auxiliary error equations and apply the discrete Sobolev embedding inequality (see Lemma \ref{lem:dis_embed}) to obtain a preliminary uniform bound result for the average value $\widetilde{\chi}_u^{n+1/2}$ (see \eqref{err_eq:inf_uv_sob}).  
Second, leveraging this bounded result, we re-estimate the auxiliary error equation system, producing a new optimal-order error analysis with a $\mK^*$-independent constant. Third, the crucial uniform bounded results for $\chi_u^{m+1}$ and $ \chi_{\phi}^{m+3/2}$ (see \eqref{err_eq:u:e1}--\eqref{err_eq:phi:e1}) are then established, by repeating the $H^2$-seminorm estimate of the auxiliary errors with the obtained new optimal-order estimates and reapplying the discrete Sobolev embedding inequality.
(ii) While for $\tau \le h^{d/2}$, the results of the $L^\infty$-norm bounds follow directly from the inverse inequality without any coupling mesh condition, as in this case ${\tau^2}/{h^{d/2}}\le \tau$ holds (see  \eqref{err_eq:u_phi:e2}). Finally, the uniform bounds of the numerical solutions imply the equivalence between the auxiliary error equation system and the original error equation system, which consequently yields the unconditional optimal-order error estimates for the LRCD scheme itself. Meanwhile, the optimal discrete $H^1$-norm error estimates for the primal variables are also proved.

The rest of the paper is organized as follows. In Section \ref{sec:LRCD}, we present the LRCD scheme for the CNLS system, establish its unique solvability, and rigorously analyze the mass and energy conservation laws and associated truncation errors in Sobolev spaces. In Section \ref{sec:AuEE}, we first derive the error estimates for the auxiliary relaxation variables at the first time level, and then construct the auxiliary error equation system via the cut-off technique. In Section \ref{sec:err_est}, we establish the optimal-order $H^1$-norm error estimates for the primal variables at the time nodes and the $L^2$-norm error estimates for the auxiliary relaxation variables at the intermediate time nodes, by firstly developing the unconditional and optimal $L^2$-norm error estimates and the $L^\infty$-bounds for the auxiliary system. In Section \ref{sec:num}, numerical examples are given to confirm the theoretical results and to demonstrate the capability of the LRCD scheme in simulating long-term wave propagation and soliton collisions. Conclusions are summarized in Section \ref{sec:con}.

\section{A linear relaxation compact difference scheme}\label{sec:LRCD}
 In this section, we propose the structure-preserving linear relaxation-based compact difference scheme for solving the CNLS system \eqref{model:e1a}--\eqref{model:e1b} on staggered temporal grids. For simplicity of presentation, we consider the description and analysis of the scheme for the CNLS system only in the two-dimensional (2D) spaces, i.e. $\Omega :=(x_{l},x_{r})\times(y_{l},y_{r})$. However, the ideas and analysis techniques used can be easily extended to the one-dimensional (1D) and three-dimensional (3D) spaces.

\subsection{Notations and preliminaries}
First, we introduce some basic notations. Let $N$ be a given positive integer and define uniform temporal grids by $t_n:=n\tau$ $(0 \le n \le N)$ with step size $\tau:=T/N$. Denote $t_{1/4}:= \tau/4$ and $t_{n+{1}/{2}} := (t_{n}+t_{n+1})/2$ for $0 \le n \le N$.
For any discrete temporal grid functions $\{w^n \mid n=0, \ldots, N\}$, we define
$$
	\widetilde{w}^{n+1/2}:=\f{w^{n+1}+w^n}{2}, \quad \delta_t w^{n+1}:=\f{w^{n+1}-w^{n}}{\tau}.
$$
    
Moreover, for given positive integers $M_x$ and $M_y$, we introduce a set of uniform spatial grids $(x_i,y_j):=(x_l+ ih_x,y_l+ jh_y)$ for $(i,j)\in \omega_h:=\{(i,j) \mid 0\le i\le M_x,~0\le j\le M_y\}$ with corresponding grid sizes $h_x:= (x_r-x_l)/M_x$ and $h_y:= (y_r-y_l)/M_y$. Let $h:=\max\{h_x,h_y\}$. Accordingly, we introduce the discrete spaces of grid functions
$$   
    \begin{aligned}
		&\mbb{C}_h:=\{v=\{v_{i,j}\}:v_{i,j}\in\mbb{C}, ~~(i,j)\in \omega_h\},~
		&\mbb{C}_h^p:=\{v\in\mbb{C}_h:v_{0,j}=v_{M_x,j}, ~v_{i,0}=v_{i,M_y}\}.\\	
		&\mbb{R}_h:=\{v=\{v_{i,j}\}:v_{i,j}\in\mbb{R},~~(i,j)\in \omega_h\},~
		&\mbb{R}_h^p:=\{v\in\mbb{R}_h:v_{0,j}=v_{M_x,j}, ~v_{i,0}=v_{i,M_y}\}.
   \end{aligned}
$$
For any  $w\in \mbb{C}_h^p$, we define the spatial difference operators and compact operator along the $x$-direction as follows:
	$$
	\delta_{x}w_{i+1/2,j} :=\f{w_{i+1,j} -w_{i,j}}{h_x},\quad\delta_{x}^{2}w_{i,j} :=\f{w_{i+1,j}-2w_{ij} +w_{i-1,j} }{h_x^{2}},
	$$
	$$
	\mA_x w_{i,j} :=	w_{i,j} +\f{h^{2}_x}{12}\delta_{x}^{2}w_{i,j} =\f{1}{12}(w_{i-1,j} +10w_{ij} +w_{i+1,j} ).
	$$
Similarly, the $y$-direction spatial difference operators and compact operator can also be defined. Furthermore, let $\mA_h:=\mA_x\mA_y$, $\Delta_h:=\delta^2_x+\delta_y^2$ and $\Lambda_h:=\mA_y\delta^2_x+\mA_x\delta_y^2.$
 Besides, for any $v,q \in \mbb{C}^p_h$ or $ \mbb{R}^p_h$,  we introduce the discrete inner products    
  $$
     \begin{aligned}
         &(v, q):= h_x h_y\sum_{i=1}^{M_x}\sum_{j=1}^{M_y} v_{i,j} \bar{q}_{i,j},~  <\delta_x v, \delta_x q>:= h_x h_y\sum_{i=1}^{M_x}\sum_{j=1}^{M_y} \delta_x v_{i-1/2,j} \delta_x \bar{q}_{i-1/2,j},  \\
         & <\delta_y v, \delta_y q>:= h_x h_y\sum_{i=1}^{M_x}\sum_{j=1}^{M_y} \delta_y v_{i,j-1/2} \delta_y \bar{q}_{i,j-1/2},
     \end{aligned}
   $$ 
and corresponding discrete norms
$$
	\|v\|:=\sqrt{(v, v)},~~	\|v\|_\mA:=\sqrt{(\mA_h v, v)},~~\|\delta_x v\|:=\sqrt{<\delta_x v, \delta_x v>}, ~~\|\delta_yv\|:=\sqrt{<\delta_y v, \delta_y v>},
$$
    $$
	|v|_1:=\sqrt{\|\delta_xv\|^2+\|\delta_yv\|^2},
	~~\|v\|_{1}:=\sqrt{\|v\|^2+|v|_1^2}, ~~\|v\|_{\infty}:=\max_{i,j} |v_{i,j}|.
	$$

Next, we introduce some useful lemmas that shall be employed in the subsequent analysis.
	The following three important lemmas can be found in Refs. \cite{LHL2021,XSS2012,ZG2023}.
\begin{lemma}\label{lem:dis_embed}
		For any $w\in\mbb{C}_h^p$ , there exists a positive constant $C_{\Omega}$ such that
		$$
		\|w \|_{\infty}\le  C_{\Omega} \left(\|w \|+\|\Lambda_h w \| \right). 
		$$
\end{lemma}
\begin{lemma}\label{lem:norm}
		For any $w\in\mbb{C}_h^p$, we have
		$$ 
		\f{4}{9}\| w\|\le\| \mA_h w\|\le\|w \|, 
           ~~\|\mA_h^{-1}w\|\le\f{9}{4}\|w\|.  
		$$
\end{lemma}
\begin{lemma}\label{lem:ieq:s2}
		For $v_a, v_b, z_a, z_b \in \mbb{C}_h$, it holds 
		\begin{equation*}\label{s2:e1}
			\begin{aligned}
				&\|  |v_a|^2-|v_b|^2 \|  \leq (\|v_a\|_{\infty}+\|v_b\|_{\infty} )\|v_a-v_b\|.
			\end{aligned}
		\end{equation*}
		Furthermore, let $\mrm{S}(v_a, v_b, z_a, z_b):=|v_a|^2-|v_b|^2-|z_a|^2+|z_b|^2$. Then, we have
		\begin{equation*}\label{s2:e2}
			\begin{aligned}
				\|\mrm{S}(v_a, v_b, z_a, z_b)\| \le  2\|z_a-z_b\|_{\infty} \|v_b-z_b\| +\mrm{G}_{\mrm{I}}(v_a, v_b, z_a, z_b) \|v_a-v_b-z_a+z_b\|,
			\end{aligned}
		\end{equation*}
		where $\mrm{G}_{\mrm{I}}(v_a, v_b, z_a, z_b):=\|v_a\|_{\infty}+\|v_b\|_{\infty}+\|z_a-z_b\|_{\infty}$.
\end{lemma}

The discrete operators $\mA_h$, $\Lambda_h$ and $\mA_h^{-1}\Lambda_h$ exhibit the following attributes.
\begin{lemma}\label{lem:operator}
		For any $w,v\in\mbb{C}_h^p$, we have
        $$
        (\mA_h w,v )=( w,\mA_h v ),~(\Lambda_h w,v )=( w,\Lambda_h v ),~(\mA_h^{-1}\Lambda_h w,v )=( w,\mal 
             A_h^{-1}\Lambda_h v ).
        $$
         Furthermore, it follows that
        $$
        (\mA_h w,w),~(\Lambda_hw,w),~(\mA_h^{-1}\Lambda_h w,w)\in \mbb{R}.
        $$
\end{lemma}
\begin{proof}
		Let $\mrm I_{\iota}$ be the $M_{\iota}$-by-$M_{\iota}$ identity matrix, and
		$$
		\mrm D_{\iota}=\f{1}{h_{\iota}^2}
		\begin{bmatrix}
			-2 & 1 & 0 & \cdots & 1 \\
			1 & -2 & 1 & \cdots & 0 \\
			\vdots & \ddots & \ddots & \ddots & \vdots \\
			1 & 0 & \cdots & 1 & -2
		\end{bmatrix},~
        \mrm A_{\iota}=\f{1}{12}
		\begin{bmatrix}
			10 & 1 & 0 & \cdots & 1 \\
			1 & 10 & 1 & \cdots & 0 \\
			\vdots & \ddots & \ddots & \ddots & \vdots \\
			1 & 0 & \cdots & 1 & 10
		\end{bmatrix}
		$$
be two $M_{\iota}$-by-$M_{\iota}$ real symmetric circulant triangular matrices (guaranteed by periodic boundary conditions) for $\iota=x, y$. Let $M:=M_xM_y$. Then, the $M$-by-$M$ matrices $\mrm A:=\mrm A_x\otimes \mrm A_y$, $\mrm B:=\mrm A_x\otimes \mrm D_y+ \mrm D_x\otimes \mrm A_y$ and $\mrm Q:=\mrm A^{-1} \mrm B$ are corresponding to the discrete operators $\mA_h$, $\Lambda_h$ and $\mA_h^{-1}\Lambda_h$, respectively. Moreover, it follows from Lemma 3.4 of \cite{FLW19} (see also \cite{HR1994}) that matrices $\mrm A$, $\mrm B$ and $\mrm Q$ are symmetric, i.e., $\mrm A = \mrm A^\top$, $\mrm B = \mrm B^\top$ and $\mrm Q = \mrm Q^\top$.
        
Below, for any $w\in\mbb{C}_h^p$, we view it as a $M$-dimensional column vector such that $w=( w_{1,1},\cdots, w_{M_x,1},\cdots,w_{1,M_y},\cdots,w_{M_x,M_y})^\top$. Then, it is easy to show that
$$
\begin{aligned}
    (\mA_h w, v ) & =h_xh_yv^{\dagger}\mrm Aw= h_xh_y [\mrm Av]^\dagger w= ( w, \mA_h v),\\
    (\Lambda_h w, v ) & =h_xh_yv^{\dagger}\mrm Bw= h_xh_y [\mrm Bv]^\dagger w= ( w, \Lambda_hv),\\
    (\mA_h^{-1}\Lambda_h w, v ) & =h_xh_yv^{\dagger}\mrm Qw= h_xh_y [\mrm Qv]^\dagger w= ( w, \mA_h^{-1}\Lambda_hv).
\end{aligned}
$$
Moreover, it further reveals that
$
(\mA_h w,w)$, $(\Lambda_hw,w)$ and $(\mA_h^{-1}\Lambda_h w,w) $ are real, which implies the second conclusion.
\end{proof}

\begin{lemma}\label{lem:H1}
          For any $w\in\mbb{C}_h^p$, we have
	    \begin{equation*}
	        \f23|w|_1^2\le-(\Lambda_hw,w)\le|w|_1^2.
	    \end{equation*}
\end{lemma}
\begin{proof}
       By using the Green formula and the inverse inequality, it follows that
       $$ (\mA_{\iota} w,w)=\|w\|^2-\f{h_{\iota}^2}{12}\|\delta_{\iota} w\|^2\ge\|w\|^2-\f13\|w\|^2=\f23\|w\|^2,$$
      which leads to $\f23\|w\|^2\le(\mA_{\iota} w,w)\le\|w\|^2$, ${\iota}=x,y.$
       Furthermore, we have
       $$\begin{aligned}
           &\f23\|\delta_xw\|^2\le-(\mA_y\delta_x^2w,w)=(\mA_y\delta_xw,\delta_xw)\le\|\delta_xw\|^2,\\
           &\f23\|\delta_yw\|^2\le-(\mA_x\delta_y^2w,w)=(\mA_x\delta_yw,\delta_yw)\le\|\delta_yw\|^2.
       \end{aligned}$$
       Thus, the conclusion follows from combinations of the above two results.
\end{proof}

 Let $\mT_h,~\mS_h:~\mbb{C}^p_h\mapsto\mbb{C}^p_h $ be two linear operators defined by $\mT_h:=\mal I_h-\f{\mrm i \kappa \tau}{2}\mA_h^{-1}\Lambda_h$, $\mS_h:=\mal I_h+\f{\mrm i \kappa \tau}{2}\mA_h^{-1}\Lambda_h$, where $\mal I_h:~\mbb{C}^p_h\mapsto\mbb{C}^p_h $ is the identity operator. It is easy to check that $\mT_h $ and $\mS_h$ are invertible. Furthermore, define $\mal C_h:=\mT_h^{-1}\mS_h$. The following two important lemmas hold; see Lemma 2.10 of \cite{ ZHOU2025} and the proof of Theorem 3.1 in \cite{ZG2023} for references.
\begin{lemma}[\cite{ZHOU2025}]\label{lem:OperatorNorm}
		For any $v\in\mbb{C}_h^p$, we have 
		$$
		\|\mal C_h v\|=\| v\|. 
		$$
\end{lemma}
\begin{lemma}[\cite{ZG2023}]\label{lem:space_operator}
   Let $\lbrace y^n \rbrace_{n=1}^N $ be a time sequence in $\mbb{C}^p_h$ satisfying
	\begin{equation*}\label{operators1}
			y^{n+1}=(\mal{C}_h-\mal{I}_h) y^n+\mal{C}_hy^{n-1}+\mT_h^{-1}\Gamma^n,\quad n=2,\ldots,N-1.
		\end{equation*}
		Then, it follows that
        $$	
          \begin{aligned}
			y^{n+1} &=\f{1}{2} \big[(-1)^{n+1}\mal{I}_h+\mal{C}_h^{n}\big]\mT_hy^2
                         +\f{1}{2}\big[(-1)^{n}\mal{C}_h+\mal{C}_h^{n}\big]\mT_hy^1
                          +\f{1}{2} \sum_{\ell=2}^n \big[(-1)^{n-\ell}\mal{I}_h+\mal{C}_h^{n+1-\ell}\big]\Gamma^\ell.
		\end{aligned}
        $$
\end{lemma}

\subsection{The fully decoupled LRCD scheme}\label{subsec:LRCD}
Introduce two auxiliary real relaxation variables $\phi=|u|^2$ and $\psi=|v|^2$. Then, the original CNLS system \eqref{model:e1a}  can be rewritten as
\begin{equation}\label{model:Relaxation}
      \left\{
		\begin{aligned}
			& \mrm{i}u_{t}+\kappa\Delta u +(\phi+\beta\psi)u=0,  &&\qquad \text{in} ~\Omega\times I,\\
			& \mrm{i}v_{t}+\kappa \Delta v+(\psi+\beta \phi)v=0,  &&\qquad \text{in} ~\Omega\times I,\\ 
			& \phi=|u|^2,	\quad \psi=|v|^2,                         &&\qquad \text{in} ~\Omega\times I. 
		\end{aligned}
      \right.
\end{equation}

Next, we shall discrete the original variables $u$ and $v$ at the time nodes $\{t_{n}\} $ and the relaxation variables $\phi$ and $\psi$ at the intermediate time nodes $\{t_{n+1/2}\}$ via the second-order Crank-Nicolson formula, combined with compact fourth-order spatial approximation of the averaged Laplacian operator $\mA_h \Delta$ by $\Lambda_h$. In particular, to maintain global second-order temporal accuracy for the two relaxation variables, a predictor-corrector scheme is developed at the first time level to predict two solutions $\{U^{1/2}, V^{1/2}\}$ for $\{\Phi^{1/2}, \Psi^{1/2}\}$. Then, the fully-discrete LRCD scheme is formulated as follows: 
\paragraph{\indent \bf Step 1} 
   Given initial values $(U^{0}, V^{0}, \Phi^{0}, \Psi^{0} ) := (u^0, v^{0}, \phi^0:=|u^{0}|^2, \psi^0:=|v^{0}|^2)$,  find $\{U^{1/2}, V^{1/2}\}\in\mbb{C}_h^p \times \mbb{C}_h^p$ such that
   \begin{align}
	&\mrm{i}\mA_h\f{U^{1/2}-U^{0}}{\tau/2}+\kappa \Lambda_h\f{U^{1/2}+U^{0}}{2}
             +\mA_h \Big[(\Phi^{0}+\beta \Psi^{0})\f{U^{1/2}+U^{0}}{2}\Big] =0,\label{LRCD:e1a} \\
	&\mrm{i}\mA_h\f{V^{1/2}-V^{0}}{\tau/2}+\kappa \Lambda_h\f{V^{1/2}+V^{0}}{2}
            +\mA_h \Big[(\Psi^{0}+\beta\Phi^{0})\f{V^{1/2}+V^{0}}{2}\Big] =0,\label{LRCD:e1b}
   \end{align}
   
\paragraph{\indent \bf Step 2}  Define $\{\Phi^{1/2}, \Psi^{1/2}\} \in \mbb{R}_h^{p}\times \mbb{R}_h^p$ by
    \begin{equation}\label{LRCD:e2a}
		\Phi^{1/2}:=|U^{1/2}|^2, \quad \Psi^{1/2}:=|V^{1/2} |^2,
    \end{equation}
and then, find $\{ U^{1}, V^{1}\} \in\mbb{C}_h^p\times \mbb{C}_h^p$ such that
	\begin{align}
		&\mrm{i}\mA_h \delta _tU^{1}+\kappa \Lambda_h\widetilde{U}^{1/2}+\mA_h \Big[(\Phi^{1/2}+\beta \Psi^{1/2})\widetilde{U}^{1/2}\Big] =0, \label{LRCD:e2b}\\
		&\mrm{i}\mA_h \delta _tV^{1}+\kappa \Lambda_h\widetilde{V}^{1/2}+\mA_h \Big[(\Psi^{1/2}+\beta\Phi^{1/2})\widetilde{V}^{1/2}\Big]=0.\label{LRCD:e2c}
	\end{align}
    
\paragraph{\indent \bf Step 3} For $n=1, \ldots, N-1$, first solve $\{\Phi^{n+1/2}, \Psi^{n+1/2}\}\in \mbb{R}^p_h \times \mbb{R}_h^p$ by
	\begin{equation}\label{LRCD:e3a}
		\Phi^{n+ 1/2}=2|U^{n}|^2-\Phi^{n- 1/2},~~ \Psi^{n+ 1/2}=2|V^{n}|^2-\Psi^{n- 1/2},
	\end{equation}
 and then,  find $\{ U^{n+1}, V^{n+1}\} \in\mbb{C}_h^p\times \mbb{C}_h^p$ such that
\begin{align}
   &\mrm{i}\mA_h \delta _tU^{n+1}+\kappa \Lambda_h\widetilde{U}^{n+ 1/2}+\mA_h \Big[(\Phi^{n+ 1/2}+\beta \Psi^{n+ 1/2})\widetilde{U}^{n+ 1/2}\Big] =0, \label{LRCD:e3b}\\
   &\mrm{i}\mA_h \delta _tV^{n+1}+\kappa \Lambda_h\widetilde{V}^{n+ 1/2}+\mA_h \Big[(\Psi^{n+ 1/2}+\beta\Phi^{n+ 1/2})\widetilde{V}^{n+ 1/2}\Big]=0.\label{LRCD:e3c}
\end{align}

\begin{remark} We remark that the above scheme possesses four essential advantages: (i) linearity in its algebraic structure;  (ii) sequential decoupling of primal and auxiliary variables, and thus, it is easy to realize parallel computing for the two primal variables at each time level; (iii) high-order accuracy in space; and  (iv) unique solvability, and mass and energy conservations at the discrete level. In practice, we first predict two  second-order accurate relaxation variables $\{\Phi^{1/2}, \Psi^{1/2}\} \in \mbb{R}_h^{p}\times \mbb{R}_h^p$  via \eqref{LRCD:e1a}--\eqref{LRCD:e2a}, and then we obtain two corrected primal solutions $\{ U^{1}, V^{1}\} \in\mbb{C}_h^p\times \mbb{C}_h^p$ via \eqref{LRCD:e2b}--\eqref{LRCD:e2c}. Next, for each $ n \ge 1$,  we first update $\{\Phi^{n+1/2}, \Psi^{n+1/2}\}$ through direct calculations of \eqref{LRCD:e3a}, thus enabling the computations of  $\{ U^{n+1}, V^{n+1}\} \in\mbb{C}_h^p\times \mbb{C}_h^p$ for the next time level via \eqref{LRCD:e3b}--\eqref{LRCD:e3c} by solving two variable-coefficient linear algebraic systems, respectively. Note that the last process can be solved in parallel, as the solutions of  $\{ U^{n+1}, V^{n+1}\}$ are decoupled with each other. 
\end{remark}

Finally, we show the existence and uniqueness of the proposed LRCD scheme.
\begin{theorem}\label{thm:unique}
    There exist unique solutions $\{U^{1/2}, V^{1/2}\}\in\mbb{C}^p_h \times \mbb{C}_h^p$ and $\{\Phi^{1/2},\Psi^{1/2}\}\in\mbb{R}^p_h \times \mbb{R}_h^p$ satisfying \eqref{LRCD:e1a}--\eqref{LRCD:e2a}, unique solutions $\{ U^{1}, V^{1}\} \in\mbb{C}_h^p\times \mbb{C}_h^p$ satisfying \eqref{LRCD:e2b}--\eqref{LRCD:e2c}, and unique solutions $\{U^{n+1}, V^{n+1}\}\in\mbb{C}^p_h \times \mbb{C}_h^p$ and $\{\Phi^{n+1/2},\Psi^{n+1/2}\}\in\mbb{R}^p_h \times \mbb{R}_h^p$ satisfying \eqref{LRCD:e3a}--\eqref{LRCD:e3c}.
\end{theorem}
\begin{proof}
     For given parameter $\nu\in\mbb{R}_h$ and constant $\varsigma>0$, let $\mal L_{\varsigma}[\nu]:\mbb{C}^p_h\mapsto\mbb{C}^p_h$ be a linear operator given by 
    \begin{equation*}
         \mal L_{\varsigma}[\nu] w:= w -\mrm i\varsigma\kappa\tau\mA_h^{-1}\Lambda_h w-\mrm i\varsigma\tau\nu w.
     \end{equation*}
     Then, it follows from Lemma \ref{lem:operator} that
      \begin{equation*}
         \mrm{Re} \big(\mal L_{\varsigma}[\nu]w,w\big)=\|w\|^2,
     \end{equation*}
     which means $\mrm{Ker}(\mal L_{\varsigma}[\nu])=\{0\}$, and thus, $\mal L_{\varsigma}[\nu]$ is invertible as $\mbb{C}^p_h$ is a finite dimensional space. 
     According to \eqref{LRCD:e1a}--\eqref{LRCD:e1b}, \eqref{LRCD:e2b}--\eqref{LRCD:e2c} and  \eqref{LRCD:e3b}--\eqref{LRCD:e3c}, we conclude that
     \begin{equation*}
        \begin{aligned}
            &U^{1/2}=(\mal L_{1/4}[\nu_1])^{-1}Z_u^0, 
              \quad V^{1/2}=(\mal L_{1/4}[\nu_2])^{-1}Z_v^0, \\
            & U^{n+1}=(\mal L_{1/2}[\nu_3])^{-1}Z_u^{n+1},
              \quad V^{n+1}=(\mal L_{1/2}[\nu_4])^{-1}Z_v^{n+1},\quad 0\le n\le N-1,
        \end{aligned}
     \end{equation*}
     with 
    \begin{equation*}
        \begin{aligned}
            &\nu_1:=\Phi^{0}+\beta \Psi^{0},~Z_u^0:=U^0+\f{\mrm i\kappa\tau}{4}\mA_h^{-1}\Lambda_h U^0+\f{\mrm i\tau}{4}\nu_1 U^0,\\
            &\nu_2:=\Psi^{0}+\beta \Phi^{0},~Z_v^0:=V^0+\f{\mrm i\kappa\tau}{4}\mA_h^{-1}\Lambda_h V^0+\f{\mrm i\tau}{4}\nu_2 V^0,\\
            &\nu_3:=\Phi^{n+1/2}+\beta \Psi^{n+1/2},~Z_u^{n+1}:=U^{n}+\f{\mrm i\kappa\tau}{2}\mA_h^{-1}\Lambda_h U^{n}+\f{\mrm i\tau}{4}\nu_1 U^{n},\\
            &\nu_4:=\Psi^{n+1/2}+\beta \Phi^{n+1/2},~Z_v^{n+1}:=V^{n}+\f{\mrm i\kappa\tau}{2}\mA_h^{-1}\Lambda_h V^{n}+\f{\mrm i\tau}{4}\nu_2 V^{n}
        \end{aligned}
     \end{equation*}
Therefore, the solutions $\{U^{1/2}, V^{1/2}\}$ and $\{U^{n+1}, V^{n+1}\}_{n \ge 0}$  are uniquely solvable. Furthermore, the explicit formulas \eqref{LRCD:e2a} and \eqref{LRCD:e3a}  directly imply the unique existence of $\{\Phi^{n+1/2},\Psi^{n+1/2}\}_{n \ge 0}$.
\end{proof}
     
\subsection{Discrete conservation laws}
In this subsection, we prove that the developed scheme \eqref{LRCD:e1a}--\eqref{LRCD:e3c} unconditionally preserves both mass and energy conservation laws at the discrete level, i.e., discrete versions of \eqref{model:conserveMass}--\eqref{model:conserveEnergy} hold.
\begin{theorem}\label{thm:conservation}
Let $\{U^{n+1},V^{n+1},\Phi^{n+1/2},\Psi^{n+1/2} \}\in \mbb{C}_h^p\times \mbb{C}_h^p\times \mathbb{R}_h^p\times \mathbb{R}_h^p$ be solutions of the LRCD scheme \eqref{LRCD:e1a}--\eqref{LRCD:e3c}. Then, the following discrete mass and energy conservation laws hold:
	\begin{equation*} 
		\mM_u^{n+1}=\mal R_u^{n+1}=\mM_u^0,~\mM_v^{n+1}=\mal R_v^{n+1}=\mM_v^0,~
		\mE^{n+1}=\mE^0, \quad \text{for}~~ n\ge 0,
	\end{equation*}
	where the discrete mass and energy are defined as 
        $$
        \begin{aligned}
            \mM_u^{n}& :=\|U^{n}\|^{2},~~\mM_v^{n}:=\|V^{n}\|^{2},~~0\le n\le N; ~~ 
            \mal{R}^{0}_u  :=\big(\Phi^{0}, 1\big),~~
            \mal{R}^{0}_v:=\big(\Psi^{0}, 1\big),\\
            \mal{R}^{n}_u & :=\f12 \big(\Phi^{n+\f12}+\Phi^{n-\f12},1\big),~~
            \mal{R}^{n}_v:=\f12 \big(\Psi^{n+\f12}+\Psi^{n-\f12},1\big),~~1\le n\le N-1;\\
            \mE^{n}&:=\left\{\begin{aligned}
 		& -\kappa\sum_{\chi=U,V}  \big(\mA_h^{-1}\Lambda_h \chi^{0}, \chi^{0}\big)
          -\f12 \left[(2\Phi^0-\Phi^{1/2},\Phi^{1/2})+ (2\Psi^0-\Psi^{1/2},\Psi^{1/2})\right]\\
 		&\quad-\f{\beta}{2} \left[(2\Phi^0-\Phi^{1/2},\Psi^{1/2})+ (2\Psi^0-\Psi^{1/2},\Phi^{1/2})\right],~~n=0, \\
 		&-\kappa \sum_{\chi=U,V}  \big(\mA_h^{-1}\Lambda_h \chi^{n}, \chi^{n}\big)-\f12\sum_{\chi=\Phi,\Psi}(\chi^{n-1/2},\chi^{n+1/2})\\
 		&\quad-\f{\beta}{2}\left[(\Phi^{n-1/2},\Psi^{n+1/2}) + (\Psi^{n-1/2},\Phi^{n+1/2})\right],~~1\le n\le N-1,\\
 		&-\kappa\sum_{\chi=U,V}  \big(\mA_h^{-1}\Lambda_h \chi^{N}, \chi^{N}\big)\\
        &\quad-\f12\left[\big(\Phi^{N-1/2},2|U^N|^2-\Phi^{N-1/2}\big)+\big(\Psi^{N-1/2},2|V^N|^2-\Psi^{N-1/2}\big)\right] \\
            &\quad -\f{\beta}{2}\left[\big(\Phi^{N-1/2},2|V^N|^2-\Psi^{N-1/2}\big)+\big(\Psi^{N-1/2},2|U^N|^2-\Phi^{N-1/2}\big)\right],~~n=N.
 	\end{aligned}\right. 
        \end{aligned}
        $$
\end{theorem}
\begin{proof} Applying the operator $\mA_h^{-1}$ to  both sides of \eqref{LRCD:e2b}--\eqref{LRCD:e2c} and \eqref{LRCD:e3b}--\eqref{LRCD:e3c} gives
	\begin{align}
		&\mrm{i}\delta _tU^{n+1}+\kappa \mA_h^{-1}\Lambda_h\widetilde{U}^{n+ 1/2}+ (\Phi^{n+ 1/2}+\beta \Psi^{n+ 1/2})\widetilde{U}^{n+ 1/2}=0, \label{LRCD:e4a}\\
		&\mrm{i} \delta _tV^{n+1}+\kappa \mA_h^{-1}\Lambda_h\widetilde{V}^{n+ 1/2}+ (\Psi^{n+ 1/2}+\beta\Phi^{n+ 1/2})\widetilde{V}^{n+ 1/2}=0,\label{LRCD:e4b}
	\end{align}
    for $0\le n \le N-1$.
	
	First, taking the discrete inner products of \eqref{LRCD:e4a} and \eqref{LRCD:e4b}  with $2 \tau\widetilde{U}^{n+1/2} $ and $2\tau \widetilde{V}^{n+1/2} $ respectively, and then choosing the corresponding imaginary parts. Note that both $\Phi^{n+ 1/2}$ and $\Psi^{n+ 1/2}$ are real, it then follows from Lemma \ref{lem:operator} that
	$$   \begin{aligned}
		&\mrm{Re} (\delta_tU^{n+1}, 2\tau \widetilde{U}^{n+1/2})= \|U^{n+1}\|^{2}-\|U^{n}\|^{2}=0,\\
		&\mrm{Re} (\delta_tV^{n+1}, 2\tau \widetilde{V}^{n+1/2})= \|V^{n+1}\|^{2}-\|V^{n}\|^{2}=0,
	\end{aligned}
        $$
	which gives the discrete mass conservation law, i.e, 
    \begin{equation}\label{eq:mass}
        \mM_u^{n+1}=\mM_u^{n}=\ldots=\mM_u^0, \quad \mM_v^{n+1}=\mM_v^{n}=\ldots=\mM_v^0.
    \end{equation} 
   Furthermore, it follows immediately from \eqref{LRCD:e3a} and \eqref{eq:mass} that
   \begin{equation}\label{eq:mass_re}
        \begin{aligned}
		&\Big(\f{\Phi ^{n+\f32}+\Phi ^{n+\f12}}2,1\Big)=\big( |U^{n+1}|^2,1\big)=\big( |U^{n}|^2,1\big)=\Big(\f{\Phi ^{n+\f12}+\Phi ^{n-\f12}}2,1\Big),\\
		&\Big(\f{\Psi^{n+\f32}+\Psi^{n+\f12}}2,1\Big)=\big( |V^{n+1}|^2,1\big)=\big( |V^{n}|^2,1\big)=\Big(\f{\Psi ^{n+\f12}+\Psi ^{n-\f12}}2,1\Big),
	\end{aligned}
   \end{equation}
    for $1\le n \le N-2$.  Therefore, we get another type discrete mass conservation law from \eqref{eq:mass_re} that
    $$
    \begin{aligned}
             \mal R^{n+1}_u=R^{n}_u=\ldots=\mal R^{0}_u=\mal M^{0}_u,~~
             \mal R^{n+1}_v=R^{n}_v=\ldots=\mal R^{0}_v=\mal M^{0}_v.
    \end{aligned}
    $$

Second, taking the discrete inner products of \eqref{LRCD:e4a} and \eqref{LRCD:e4b} with $U^{n+1}-U^{n} $ and $V^{n+1}-V^{n} $ respectively, and then choosing the corresponding real parts and summing them together yields
\begin{equation}\label{LRCD:e5}
   \begin{aligned}
	& -\kappa \sum_{\chi=U,V}  \left[\big(\mA_h^{-1}\Lambda_h \chi^{n+1}, \chi^{n+1}\big) -  \big(\mA_h^{-1}\Lambda_h \chi^{n}, \chi^{n}\big)\right] \\
	& ~~ =  \big(\Phi^{n+1/2} + \beta \Psi^{n+1/2}, |U^{n+1}|^2-|U^{n}|^2\big)
	     + \big(\Psi^{n+1/2} + \beta \Phi^{n+1/2}, |V^{n+1}|^2-|V^{n}|^2\big),
	\end{aligned}
\end{equation}
    for $0\le n \le N-1$. As for $1\le n\le N-2$, note that it follows immediately from \eqref{LRCD:e3a} that
 $$
   |U^{n+1}|^2-|U^{n}|^2=\f{\Phi^{n+3/2}-\Phi^{n-1/2}}{2},~~  |V^{n+1}|^2-|V^{n}|^2=\f{\Psi^{n+3/2}-\Psi^{n-1/2}}{2}.
 $$
 Thus, by inserting the above formulas into the right-hand side of \eqref{LRCD:e5}, we have
  $$\begin{aligned}
 	&(\Phi^{n+1/2} + \beta \Psi^{n+1/2}, |U^{n+1}|^2-|U^{n}|^2)  + (\Psi^{n+1/2} + \beta \Phi^{n+1/2}, |V^{n+1}|^2-|V^{n}|^2)\\
        &=\f12\sum_{\chi=\Phi,\Psi}(\chi^{n+1/2},\chi^{n+3/2})+\f{\beta}{2}\left[(\Phi^{n+1/2},\Psi^{n+3/2})+(\Psi^{n+1/2},\Phi^{n+3/2})\right]\\
        &\quad-\f12\sum_{\chi=\Phi,\Psi}(\chi^{n-1/2},\chi^{n+1/2})-\f{\beta}{2}\left[(\Phi^{n-1/2},\Psi^{n+1/2})+(\Psi^{n-1/2},\Phi^{n+1/2})\right],
 \end{aligned}$$
 which together with \eqref{LRCD:e5} gives $\mE^{N-1}=\mE^{N-2}=\cdots=\mE^1$. Moreover, inserting \eqref{LRCD:e3a} with  $n =1$ into \eqref{LRCD:e5} with  $n =0$, we have 
 $$
   \begin{aligned}
       &-\kappa\sum_{\chi=U,V}  \big(\mA_h^{-1}\Lambda_h \chi^{1}, \chi^{1}\big) -\f12\sum_{\chi=\Phi,\Psi}(\chi^{1/2},\chi^{3/2})-\f{\beta}{2}\left[(\Phi^{1/2},\Psi^{3/2})+(\Psi^{1/2},\Phi^{3/2})\right]\\
       &=-\kappa\sum_{\chi=U,V}  \big(\mA_h^{-1}\Lambda_h \chi^{0}, \chi^{0}\big)-\f12\left[(2\Phi^0-\Phi^{1/2},\Phi^{1/2})+(2\Psi^0-\Psi^{1/2},\Psi^{1/2})\right]\\
       &\quad-\f{\beta}{2}\left[(2\Phi^0-\Phi^{1/2},\Psi^{1/2})+(2\Psi^0-\Psi^{1/2},\Phi^{1/2})\right],
   \end{aligned}
 $$
 which leads to $\mE^1=\mE^0$. Finally, for $n=N-1$, inserting \eqref{LRCD:e3a} into \eqref{LRCD:e5} gives
 $$\begin{aligned}
            &-\kappa\sum_{\chi=U,V}  \big(\mA_h^{-1}\Lambda_h \chi^{N-1}, \chi^{N-1}\big)
            \\
            &\quad-\f12\sum_{\chi=\Phi,\Psi}(\chi^{N-3/2},\chi^{N-1/2})
            -\f{\beta}{2}\left[(\Phi^{N-3/2},\Psi^{N-1/2})+(\Psi^{N-3/2},\Phi^{N-1/2})\right]\\
            &=-\kappa\sum_{\chi=U,V}  \big(\mA_h^{-1}\Lambda_h \chi^{N}, \chi^{N}\big)\\
        &\quad-\f12\left[\big(\Phi^{N-1/2},2|U^N|^2-\Phi^{N-1/2}\big)+\big(\Psi^{N-1/2},2|V^N|^2-\Psi^{N-1/2}\big)\right] \\
            &\quad -\f{\beta}{2}\left[\big(\Phi^{N-1/2},2|V^N|^2-\Psi^{N-1/2}\big)+\big(\Psi^{N-1/2},2|U^N|^2-\Phi^{N-1/2}\big)\right],
 \end{aligned}$$
which implies $\mE^{N-1}=\mE^N$.  Therefore, the discrete energy conservation law also holds.
\end{proof}

\begin{remark}
It is worth noting that  if we directly set $\{\Phi^{\f12},\Psi^{\f12}\}=\{\Phi^{0},\Psi^{0}\}=\{|u^{0}|^2, |v^{0}|^2\}$ in \eqref{LRCD:e2a} without a second-order accurate predicition in \eqref{LRCD:e1a}--\eqref{LRCD:e1b}, the initial discrete energy $\mE^0$ reduces to 
    $$\begin{aligned}
         \mE^{0}&:=-\kappa\left[\left(\mA_h^{-1}\Lambda_h U^{0}, U^{0}\right) +  \left(\mA_h^{-1}\Lambda_h V^{0}, V^{0}\right)\right] \\
               &\qquad-\f12 \left[ (|U^{0}|^2,|U^{0}|^2)+ (|V^{0}|^2,|V^{0}|^2)\right] -\beta \left( |U^{0}|^2, |V^{0}|^2\right),
    \end{aligned}$$
    which is consistent with the discrete original energy \eqref{model:conserveEnergy} at $t=t_0$. However, this choice  yields a second order convergence of the numerical method at the time nodes but a first order convergence at the intermediate  time nodes (see \cite{ZG2023}). 
\end{remark}

\subsection{Truncation errors}
Throughout this paper, we make the following regularity assumptions for the solutions to the CNLS model \eqref{model:Relaxation} that
\begin{equation}\label{assu:Regular}
\begin{aligned}
     &u,v\in  H^{4}(0,T;L^{2}(\Omega))\cap H^{3}(0,T;H^{2}(\Omega))\cap W^{1,\infty}(0,T;H^{6}(\Omega)).
\end{aligned}
\end{equation}
Furthermore, we assume that there is a positive constant $\mK$ such that 
\begin{equation}\label{bound_assu}
    \max_{\chi=u,v,\phi,\psi}\{\|\chi\|_{L^\infty(\Omega\times I)},\|\p_t\chi\|_{L^\infty(\Omega\times I)}\}\le \mK.
\end{equation}

Let $\{u^n,~v^n\}\in \mbb{C}_h^p$ and $\{\phi^{n+1/2}, \psi^{n+1/2}\}\in \mathbb{R}_h^p $ be the exact solutions at $t=t_{n}$ and $t=t_{n+1/2}$, respectively. Then, applying the compact operator $\mA_h$ to the first two equations in \eqref{model:Relaxation}, it is easy to verify that the following consistency equations hold at corresponding temporal nodes: 
\begin{equation}\label{model:C1}
     \mrm{i}\mA_h \f{u^{1/2}-u^0}{\tau/2}+\kappa\Lambda_h\f{u^{1/2}+u^0}{2}
        +\mA_h \Big[(\phi^{0}+\beta\psi^{0})\f{u^{1/2}+u^0}{2}\Big]+\mrm{r}^{1/4}_u=0, 
\end{equation}
\begin{equation}\label{model:C2}
     \mrm{i}\mA_h \f{v^{1/2}-v^0}{\tau/2}+\kappa\Lambda_h\f{v^{1/2}+v^0}{2}
        +\mA_h \Big[(\psi^{0}+\beta\phi^{0})\f{v^{1/2}+v^0}{2}\Big]+\mrm{r}^{1/4}_v=0,
\end{equation}
at $t=t_{1/4}$, and
\begin{equation}\label{model:C3}
    \phi^{1/2}=|u^{1/2}|^2,~~\psi^{1/2}=|v^{1/2}|^2,
\end{equation}
\begin{equation}\label{model:C4}
    \mrm{i}\mA_h\delta_tu^{1}+\kappa\Lambda_h\widetilde{u}^{1/2}
       +\mA_h\Big[(\phi^{1/2}+\beta\psi^{1/2})\widetilde{u}^{1/2}\Big]+\mrm{r}^{1/2}_u=0,
\end{equation}
\begin{equation}\label{model:C5}
   \mrm{i}\mA_h\delta_t v^{1}+\kappa\Lambda_h\widetilde{v}^{1/2}+\mA_h\Big[(\psi^{1/2}+\beta\phi^{1/2})\widetilde{v}^{1/2}\Big]+\mrm{r}^{1/2}_v=0, 
\end{equation}
at $t=t_{1/2}$, and
\begin{equation}\label{model:C6}
\phi^{n+1/2}=2|u^{n}|^2-\phi^{n-1/2}+\mrm{r}^{n}_\phi, ~~ 
    \psi^{n+1/2}=2|v^{n}|^2-\psi^{n-1/2}+\mrm{r}^{n}_\psi,  
\end{equation}
at $t=t_{n}$, and
\begin{equation}\label{model:C7}
     \mrm{i}\mA_h\delta_tu^{n+1}+\kappa\Lambda_h\widetilde{u}^{n+1/2}
       +\mA_h\Big[(\phi^{n+1/2}+\beta\psi^{n+1/2})\widetilde{u}^{n+1/2}\Big]+\mrm{r}^{n+1/2}_u=0,
\end{equation}
\begin{equation}\label{model:C8}
  \mrm{i}\mA_h\delta_t v^{n+1}+\kappa\Lambda_h\widetilde{v}^{n+1/2}+\mA_h\Big[(\psi^{n+1/2}+\beta\phi^{n+1/2})\widetilde{v}^{n+1/2}\Big]+\mrm{r}^{n+1/2}_v=0, 
\end{equation}
at $t=t_{n+1/2}$, $1\le n\le N-1$, where $\{\mrm{r}^{1/4}_u,~\mrm{r}^{1/4}_v\},~
\{\mrm{r}^{1/2}_u,~\mrm{r}^{1/2}_v\},~\{\mrm{r}^{n}_\phi,~\mrm{r}^{n}_\psi\},~\{\mrm{r}^{n+1/2}_u,~\mrm{r}^{n+1/2}_v\},$ represent the temporal and spatial truncation errors at $t=t_{1/4}$, $t=t_{1/2}$, $t=t_{n}$, and $t=t_{n+1/2}$, respectively, and 
 \begin{equation*}  
   \begin{aligned}    
        &
	\mrm{r}^{1/4}_u:= \mrm{i}\mA_h \Big( u_t^{1/4}- \f{u^{1/2}-u^0}{\tau/2}\Big)+\kappa\Big(\mA_h\Delta u^{1/4}-\Lambda_h\f{u^{1/2}+u^0}{2}\Big)\\
	&\quad\qquad+\mA_h\Big[(\phi^{1/4}+\beta\psi^{1/4})u^{1/4} -(\phi^{0}+\beta\psi^{0})\f{u^{1/2}+u^0}{2}\Big],\\
	&
	\mrm{r}^{1/4}_v:=  \mrm{i}\mA_h \Big( v_t^{1/4}- \f{v^{1/2}-v^0}{\tau/2}\Big)+\kappa\Big(\mA_h\Delta v^{1/4}-\Lambda_h\f{v^{1/2}+v^0}{2}\Big)\\
	&\quad\qquad +\mA_h\Big[(\psi^{1/4}+\beta\phi^{1/4})v^{1/4}-(\beta\phi^{0}+\psi^{0})\f{v^{1/2}+v^0}{2}\Big],\\
	&
	\mrm{r}^{n+1/2}_u:= \mrm{i}\mA_h \Big(u_t^{n+1/2}- \delta_tu^{n+1}\Big)+\kappa\Big(\mA_h\Delta u^{n+1/2}- \Lambda_h\widetilde{u}^{n+1/2}\Big)\\
       &\qquad\qquad  +\mA_h\Big[(\phi^{n+1/2}+\beta\psi^{n+1/2})({u}^{n+1/2}-\widetilde{u}^{n+1/2})\Big],~0\le n\le N-1,\\
	&
	\mrm{r}^{n+1/2}_v:= \mrm{i}\mA_h  \Big(v_t^{n+1/2}- \delta_tv^{n+1}\Big)+\kappa\Big(\mA_h\Delta v^{n+1/2}-\Lambda_h\widetilde{v}^{n+1/2}\Big)\\
      &\qquad\qquad +\mA_h\Big[(\psi^{n+1/2}+\beta\phi^{n+1/2})({v}^{n+1/2}-\widetilde{v}^{n+1/2})\Big],~0\le n\le N-1,\\
	&\mrm{r}^{n}_\phi:=\phi^{n+1/2}+\phi^{n-1/2}-2\phi^{n},~~\mrm{r}^{n}_\psi:= \psi^{n+1/2}+\psi^{n-1/2}-2\psi^{n},~1\le n\le N-1.
\end{aligned}
    \end{equation*} 
        
The following theorem establishes the orders of truncation errors.
\begin{lemma}\label{lem:trunc_errors}
Assume that the solutions to \eqref{model:e1a}--\eqref{model:e1b} satisfy the regularity condition \eqref{assu:Regular}. Then, there exists a positive constant $\widehat{C}_0$, related to the corresponding norms of the exact solutions, but independent of $\tau$ and $h$, such that the truncation errors of the LRCD scheme satisfy
\begin{equation*}
     \begin{aligned}
         &\quad\|\mrm{r}^{1/4}_u\|+\|\mrm{r}^{1/4}_v\|\le \widehat{C}_0 (\tau+h^4),\\
          &\max_{0\le n\le N-1}\{\|\mrm{r}^{n+1/2}_u \|+\|\mrm{r}^{n+1/2}_v \|\}\le \widehat{C}_0 (\tau^2+h^4),\\
          &\max_{1\le n\le N-1}\{\|\mrm{r}^{n}_\phi \|+ \|\mrm{r}^{n}_\psi \|+\|\mrm{r}^{n}_\phi \|_\infty+ \|\mrm{r}^{n}_\psi \|_\infty\}\le \widehat{C}_0 \tau^2,\\
          &\max_{1\le n\le N-1}\{(\| \mrm{r}^{n+1/2}_u- \mrm{r}^{n-1/2}_u\|+\| \mrm{r}^{n+1/2}_v- \mrm{r}^{n-1/2}_v\|\}\le \widehat{C}_0 \tau(\tau^2+h^4),\\
         &\max_{2\le n\le N-1}\{\|\mrm{r}_{\phi}^{n}-\mrm{r}_{\phi}^{n-1} \|+\|\mrm{r}_{\psi}^{n}-\mrm{r}_{\psi}^{n-1}\}\le \widehat{C}_0 \tau^3.
      \end{aligned}  
 \end{equation*}
\end{lemma}
\begin{proof} It follows directly from Lemma \ref{lem:norm} that
$$   
\begin{aligned}
   \|\mrm{r}^{1/4}_u \|
		&\le \big\| u_t^{1/4}- \f{u^{1/2}-u^0}{\tau/2} \big\| 
                 + \kappa \big\| \Delta (u^{1/4}-\f{u^{1/2}+u^0}{2}) \big\|\\
		&\quad + \kappa \big\|(\mA_h \Delta-\Lambda_h)\f{u^{1/2}+u^0}{2}  \big\|
                    + \big\|(\phi^{1/4}+\beta\psi^{1/4}-\phi^{0}-\beta\psi^{0})u^{1/4}  \big\|\\
		&\quad + \big\|(\phi^{0}+\beta\psi^{0})( \f{u^{1/2}+u^0}{2}- u^{1/4})   \big\| 
		=: \sum_{i=1}^5 \|J^{1/4}_{u,i}\|. 
\end{aligned} 
$$
Moreover, under the regularity condition, by using the Taylor expansion, it is easy to see
$$
	 \|J^{1/4}_{u,1}\| \le C\tau^2,~ \|J^{1/4}_{u,2}\| \le  C\tau^2,~ \|J^{1/4}_{u,4}\| \le C\tau,~ 
     \|J^{1/4}_{u,5}\| \le C\tau^2.
$$
Meanwhile, it is seen from Section 3.1 of \cite{ZBY2024} that
        $
           \|J^{1/4}_{u,3}\| \le Ch^4.
        $
        Therefore, we arrive at
        $$
	 \|\mrm{r}^{1/4}_u \|\le \widehat{C}_0 (\tau+h^4). 
	$$
   
   Analogously, we can show the other conclusions; see also subsection 3.3 of \cite{ZHOU2025} for reference.
Thus, the lemma is proved.
\end{proof}
\section{Auxiliary error equations}\label{sec:AuEE}
In this section, we first present the error equations of the LRCD scheme, and in particular, we discuss the error estimates for the prediction step \eqref{LRCD:e1a}--\eqref{LRCD:e2a}. Then, by introducing two novel auxiliary cut-off coefficients, we construct some auxiliary error equations with provable unique solvability, which also provide the theoretical foundation for both the uniform boundedness and unconditional optimal-order error estimates of the LRCD scheme.

Define the error functions
\begin{equation*}\label{def:err1}
 \begin{aligned}
    e_u^{1/2}:=U^{1/2}-u^{1/2}, ~ e^n_u:=U^{n}-u^n, ~e^{n+1/2}_\phi:=\Phi^{n+1/2}-\phi^{n+1/2},\\ 
    e_v^{1/2}:=V^{1/2}-v^{1/2}, ~e^n_v:=V^{n}-v^n, ~e^{n+1/2}_\psi:=\Psi^{n+1/2}-\psi^{n+1/2}.
 \end{aligned}
\end{equation*}
There hold $e^0_u=0$ and $e^0_v=0$. Moreover, by subtracting \eqref{model:C1}--\eqref{model:C8}  from \eqref{LRCD:e1a}--\eqref{LRCD:e3c} and applying the operator $\mA_h^{-1}$ on the both sides of the results, we get the following error equations:
  \begin{equation}
        \mrm i e_u^{1/2}+\f{\kappa\tau}{4}\mA_h^{-1}\Lambda_he_u^{1/2}+\f{\tau}{4}\Big[(\Phi^0+\beta\Psi^0)e_u^{1/2}\Big]
         =\f{\tau}{2}\mA_h^{-1}\mrm r^{1/4}_u,\label{err_eq:eu0}
    \end{equation}
    \begin{equation}
        \mrm i e_v^{1/2}+\f{\kappa\tau}{4}\mA_h^{-1}\Lambda_he_v^{1/2}+\f{\tau}{4}\Big[(\beta\Phi^0+\Psi^0)e_v^{1/2}\Big]
        =\f{\tau}{2}\mA_h^{-1}\mrm r^{1/4}_v,\label{err_eq:ev0}
    \end{equation}
\begin{equation}\label{err_eq:ephi0}    
  e^{1/2}_\phi   = | u^{1/2}|^2-| U^{1/2}|^2, ~~e^{1/2}_\psi   = | v^{1/2}|^2-| V^{1/2}|^2,
  \end{equation}
  and
\begin{equation}
    \mrm i\delta_{t}e^{n+1}_u+\kappa\mA_h^{-1}\Lambda_h\widetilde{e}_u^{n+1/2}+\Gamma_{u,1}^{n+1}+\Gamma_{u,2}^{n+1}+\Gamma_{u,3}^{n+1}=0,
	\label{err_eq:e1}
\end{equation}
\begin{equation}
   \mrm i\delta_{t}e^{n+1}_v+\kappa\mA_h^{-1}\Lambda_h\widetilde{e}_v^{n+1/2}+\Gamma_{v,1}^{n+1}+\Gamma_{v,2}^{n+1}+\Gamma_{v,3}^{n+1}=0,
	\label{err_eq:e2}
\end{equation}
   for $0\le n \le N-1$, and
\begin{equation}
   e^{n+1/2}_{\phi}+e^{n-1/2}_{\phi}=2(|u^{n} |^2-|U^{n}|^2  )-\mrm{r}_{\phi}^{n},
	\label{err_eq:e3}
\end{equation}
\begin{equation}
   e^{n+1/2}_{\psi}+e^{n-1/2}_{\psi}=2(|v^{n} |^2-|V^{n} |^2  )-\mrm{r}_{\psi}^{n},
	\label{err_eq:e4}
\end{equation}
 for $1\le n \le N-1$,  where, 
$$   
\begin{aligned}
	&\Gamma_{u,1}^{n+1}:= (e_\phi^{n+1/2}+\beta e_\psi^{n+1/2} )\widetilde{U}^{n+1/2},~
        &\Gamma_{v,1}^{n+1}&:= (e_\psi^{n+1/2}+\beta e_\phi^{n+1/2} )\widetilde{V}^{n+1/2}, \\
	&\Gamma_{u,2}^{n+1}:= ( \phi^{n+1/2}+\beta\psi^{n+1/2}) \widetilde{e}^{n+1/2}_u,~
        &\Gamma_{v,2}^{n+1}&:= ( \psi^{n+1/2}+\beta\phi^{n+1/2}) \widetilde{e}^{n+1/2}_v,\\
	&\Gamma_{u,3}^{n+1}:=-\mA_h^{-1}\mrm{r}_{u}^{n+1/2},~
        &\Gamma_{v,3}^{n+1}&:=-\mA_h^{-1}\mrm{r}_{v}^{n+1/2}.
\end{aligned} 
$$
From Theorem \ref{thm:unique}, we conclude that \eqref{err_eq:eu0}--\eqref{err_eq:e4} are uniquely solvable.
Next, we shall provide optimal-order error estimates for $\|e^{1/2}_\phi \|$ and $\|e^{1/2}_\psi \|$ in \eqref{err_eq:ephi0}.
\begin{lemma}\label{lem:err_ana_pre} Let $\{U^{1/2}, V^{1/2},\Phi^{1/2},\Psi^{1/2}\}\in\mbb{C}^p_h \times \mbb{C}^p_h \times \mbb{R}^p_h \times \mbb{R}^p_h$ be solutions of the prediction step \eqref{LRCD:e1a}--\eqref{LRCD:e2a}  of the LRCD scheme. There is a positive constant $\widehat{C}_1$, independent of $\tau$ and $h$, such that
\begin{equation*}\label{errphi1/2}
		\|e^{1/2}_\phi \|+\|e^{1/2}_\psi \|\le \widehat{C}_1 (\tau^2+h^4),~~\|e^{1/2}_\phi \|_\infty+\|e^{1/2}_\psi \|_\infty\le \widehat{C}_1  (\tau+h^3),
\end{equation*}	
   where  $\widehat{C}_1$ is related to the constants $\widehat{C}_0$, $\mK$ and $C_{\Omega}$.    
\end{lemma}
\begin{proof} On one hand, taking inner products of \eqref{err_eq:eu0} with $e_u^{1/2}$ and \eqref{err_eq:ev0} with $e_v^{1/2}$, then  keeping the imaginary parts of the resulting equations, we have
$$  
    \begin{aligned}
		\|e^{1/2}_u\|^2=\f{\tau}{2}\mrm{Im}(\mA_h^{-1} \mrm r^{1/4}_u,e_u^{1/2}),~~
		\|e^{1/2}_v\|^2=\f{\tau}{2}\mrm{Im}( \mA_h^{-1}\mrm r^{1/4}_v,e_v^{1/2}).
    \end{aligned} 
$$
It then follows from Cauchy-Schwarz inequality and Lemmas \ref{lem:norm}, \ref{lem:trunc_errors} that
    \begin{equation}\label{err:uv12}
        \|e^{1/2}_u \| \le \f{9\tau}{8}\| \mrm{r}^{1/4}_u\| \le \widehat{C}_0 \tau ( \tau+h^4), ~~ 
        \|e^{1/2}_v \| \le \f{9\tau}{8}\| \mrm{r}^{1/4}_v\| \le \widehat{C}_0 \tau ( \tau+h^4).
    \end{equation}

On the other hand, by taking norm on the both sides of \eqref{err_eq:eu0} and \eqref{err_eq:ev0}, and using Lemmas \ref{lem:norm}, \ref{lem:trunc_errors} and \eqref{err:uv12} again, we see
    \begin{equation}\label{err:eu12h}
        \|\Lambda_he_u^{1/2}\|\le C (\tau^{-1}\|e^{1/2}_u \|+\|e^{1/2}_u \|+\|\mrm r^{1/4}_u \|)
         \le \widehat{C}_1  ( \tau+h^4),~~
        \|\Lambda_he_v^{1/2}\| \le  \widehat{C}_1  ( \tau+h^4).
    \end{equation}

Now, combining \eqref{err:uv12} and \eqref{err:eu12h} together, and using Lemma \ref{lem:dis_embed} we have
  \begin{equation}\label{err:eu12t}
      \|e^{1/2}_u\|_{\infty}\le C_{\Omega}\big(\|e^{1/2}_u\|+ \|\Lambda_he_u^{1/2}\|\big)\le \widehat{C}_1  ( \tau+h^4), ~~
      \|e^{1/2}_v\|_{\infty}\le  \widehat{C}_1  ( \tau+h^4).
  \end{equation}
Therefore, for sufficiently small $\tau$ and $h$, we conclude from \eqref{bound_assu} and \eqref{err:eu12t} that
\begin{equation}\label{err:bound}
    \|U^{1/2} \|_\infty\le \mK+1 ,~~  \|V^{1/2} \|_\infty\le \mK+1.
\end{equation}
Furthermore, we apply Lemma \ref{lem:ieq:s2}, \eqref{err:uv12} and \eqref{err:bound} to \eqref{err_eq:ephi0} to derive
$$ \begin{aligned}
   \|e^{1/2}_\phi \|  	\le   (\|u^{1/2} \|_\infty+\|U^{1/2} \|_\infty )\|e^{1/2}_u \|	\le \widehat{C}_1  (\tau^2+h^4) ~~\text{and}~~ 
   \|e^{1/2}_\psi  \| \le \widehat{C}_1 (\tau^2+h^4).
	\end{aligned} $$
Similarly, we can obtain
$$ \begin{aligned}
		\|e^{1/2}_\phi \|_\infty\le   (\|u^{1/2} \|_\infty+\|U^{1/2} \|_\infty )\|e^{1/2}_u \|_\infty\le \widehat{C}_1(\tau+h^4) ~\text{and}~ 
		\|e^{1/2}_\psi  \|_\infty \le\widehat{C}_1(\tau+h^4).
\end{aligned} $$
Thus, the lemma is proved.
\end{proof}

 To establish the \textit{unconditional} optimal-order error estimates for the proposed LRCD scheme, and meanwhile, to avoid the difficulties arising from the requirement of $L^\infty$-norm bounds of the numerical solutions, we shall introduce an auxiliary error equation system by using four specially defined cut-off coefficients.

Let $\epsilon$ be a given small positive constant such that $0<\epsilon \ll 1$, and $\mK^*$ be a to-be-determined positive constant. Let
$\{\chi_u^{n+1}, \chi_v^{n+1},\chi_\phi^{n+1/2},\chi_\psi^{n+1/2}\}\in \mbb{C}_h^p\times \mbb{C}_h^p\times \mathbb{R}_h^p\times \mathbb{R}_h^p$ be solutions to the following auxiliary error equations that are analogous to \eqref{err_eq:e1}--\eqref{err_eq:e4}:
\begin{equation}
\mrm i\delta_{t}\chi^{n+1}_u+\kappa\mA^{-1}_h\Lambda_h\widetilde{\chi}_u^{n+1/2}+\mrm D_{u,1}^{n+1}+\mrm D_{u,2}^{n+1}+\mrm D_{u,3}^{n+1}=0,\label{err_eq:M1}
\end{equation}
\begin{equation}
\mrm i \delta_{t}\chi^{n+1}_v+\kappa\mA^{-1}_h\Lambda_h\widetilde{\chi}_v^{n+1/2}+\mrm D_{v,1}^{n+1}+\mrm D_{v,2}^{n+1}+\mrm D_{v,3}^{n+1}=0, \label{err_eq:M2}
\end{equation}
for $0\le n \le N-1$, and
\begin{equation}
    \chi^{n+1/2}_{\phi}+\chi^{n-1/2}_{\phi}=2(|u^{n} |^2-|U_\mfm^{n} |^2  )-\mrm{r}_{\phi}^{n},\label{err_eq:M3}
\end{equation}
\begin{equation}
\chi^{n+1/2}_{\psi}+\chi^{n-1/2}_{\psi}=2(|v^{n} |^2-|V_{\mfm}^{n} |^2  )-\mrm{r}_{\psi}^{n},\label{err_eq:M4}
\end{equation}
for $1\le n \le N-1$, with $\chi_u^0=e_u^0=0$, $\chi_v^0=e_v^0=0$, $\chi_\phi^{1/2}=e_\phi^{1/2}$, $\chi_\psi^{1/2}=e_\psi^{1/2}$,  and
$$   \begin{aligned}
	&\mrm D_{u,1}^{n+1}= (\mfm_\phi\chi^{n+1/2}_\phi+\beta \mfm_\psi\chi^{n+1/2}_\psi )\widetilde{U}^{n+1/2}_{{\mfm}},
	&&\mrm D_{v,1}^{n+1}=(\mfm_\psi\chi^{n+1/2}_\psi+\beta \mfm_\phi\chi^{n+1/2}_\phi )\widetilde{V}^{n+1/2}_{{\mfm}}, \\
	&\mrm D_{u,2}^{n+1}=( \phi^{n+1/2}+\beta\psi^{n+1/2})\, \mfm_u\widetilde{\chi}_u^{n+1/2} ,
	&&\mrm D_{v,2}^{n+1}=( \psi^{n+1/2}+\beta\phi^{n+1/2})\,\mfm_v\widetilde{\chi}_v^{n+1/2}, \\
	&\mrm D_{u,3}^{n+1}=-\mA^{-1}_h\mrm{r}_{u}^{n+1/2}, &&\mrm D_{v,3}^{n+1}=-\mA^{-1}_h\mrm{r}_{v}^{n+1/2}.
\end{aligned} $$ 
Here, the cut-off coefficients 
 \begin{equation}\label{cut_off:e1}
    \begin{aligned}
     \mfm_u:=\min\Big\{\f{1}{\max_n\|{\chi}_u^{n}\|_\infty+\epsilon},1\Big\},~~\mfm_\phi:=\min\Big\{\f{\mK^*}{\max_n\|{\chi}_\phi^{n+1/2}\|_\infty+\epsilon},1\Big\},\\
     \mfm_v:=\min\Big\{\f{1}{\max_n\|{\chi}_v^{n}\|_\infty+\epsilon},1\Big\}, ~~\mfm_\psi:=\min\Big\{\f{\mK^*}{\max_n\|{\chi}_\psi^{n+1/2}\|_\infty+\epsilon},1\Big\},
    \end{aligned} 
 \end{equation}
 and the auxiliary variables
\begin{equation}\label{cut_off:e2}
    \begin{aligned}
       U^n_\mfm:=u^n+\mfm_u\chi^n_u,~~\Phi^{n+1/2}_\mfm:=\phi^{n+1/2}+\mfm_\phi\chi^{n+1/2}_\phi, \\
       V^n_\mfm:=v^n+\mfm_v\chi^n_v,~~\Psi^{n+1/2}_\mfm:=\psi^{n+1/2}+\mfm_\psi\chi^{n+1/2}_\psi.
    \end{aligned} 
\end{equation}

\begin{remark} As the uniform boundedness of the numerical solutions cannot be assured unconditionally, it is not easy to prove the unconditional optimal-order error estimates for the LRCD scheme \eqref{err_eq:e1}--\eqref{err_eq:e4} directly in multiple space dimensions \cite{ZHOU2025,ZG2023}. In fact, for sufficiently small $\tau$ and $h$, we shall show that  $\mfm_u=\mfm_v =\mfm_\phi=\mfm_\psi \equiv 1$ in \eqref{cut_off:e1}. In this context, the auxiliary error equations \eqref{err_eq:M1}--\eqref{err_eq:M4} are completely equivalent to the original \eqref{err_eq:e1}--\eqref{err_eq:e4}, and therefore, 
they have the same uniquely solvable numerical solutions. Thus, instead of analyzing the original LRCD scheme itself, we shall discuss the error estimates for the auxiliary error equations \eqref{err_eq:M1}--\eqref{err_eq:M4} in the following. The similar idea can also be found in, for example, Refs. \cite{XXF'22,ZBY2024,ZHOU2025,ZG2023}.
\end{remark}

\begin{remark} Based on the definitions \eqref{cut_off:e1}--\eqref{cut_off:e2}, the boundedness results for the auxiliary error equations \eqref{err_eq:M1}--\eqref{err_eq:M4} can be directly implied that
\begin{equation}\label{ieq:bound_UmVm}
    \|U^n_\mfm \|_\infty\le \mathcal{K}+1,~~   \|V^n_\mfm \|_\infty\le \mathcal{K}+1,
 \end{equation}
 and
  \begin{equation}\label{ieq:bound_phi}
  \|\mfm_\phi\chi_\phi^{n+1/2} \|_\infty\le \mK^*,~~    \|\mfm_\psi\chi_\psi^{n+1/2} \|_\infty\le \mK^*,
 \end{equation}
 which shall paly an important role in the unconditional optimal-order error analysis.
\end{remark}

\section{Unconditional optimal-order error estimates} \label{sec:err_est}
In this section, we aim to show the unconditional optimal-order error estimates for the LRCD scheme, relying on optimal-order error estimates for the auxiliary error equations \eqref{err_eq:M1}--\eqref{err_eq:M4} and the equivalence between \eqref{err_eq:M1}--\eqref{err_eq:M4} and \eqref{err_eq:e1}--\eqref{err_eq:e4}. 
To this aim, we introduce some specific positive constants that independent of mesh parameters as follows:
 \begin{equation}\label{p:constant}
   C_1 := 2(1+|\beta|)(\mK+1),~  C_2 := 16\mK + 8,~  C_3 := \f{C_1 C_2}{4},~ C_4 :=3C_1+2C_3,
\end{equation}
\begin{equation}\label{p:constant-star}
       C_1^* := \f{(1+|\beta|)}{2}\mK^*,~   C_2^* := 2C_1^* + C_3,~   C_3^* :=2\max\{2C_1+C_3,3C_2^*+3C_1\}.
\end{equation}
Note that the constants in \eqref{p:constant-star} depend on the to-be-determined constant $\mathcal{K}^*$. Besides, we also define some generic positive constants $\{\widehat{C}_i\}$ ($0 \le i \le 5$) that are independent of both mesh parameters and $\mathcal{K}^*$,  but may rely on the bound $\mK$, the final time $T$ and the exact solutions in corresponding norms. 


The subsequent lemma presents error estimates of $ \|\chi^{m+1}_u \|$ and $\|\chi^{m+1}_v \|$, which are shown to be related to those of $\|\chi_\phi^{n+1/2} \|$ and $\|\chi_\psi^{n+1/2} \|$ for $0 \le n \le m$.
\begin{lemma}\label{lem:M_err1} 
    Under the regularity condition \eqref{assu:Regular}, there exists a positive constant $\widehat{C}_2$ related to $\widehat{C}_0$ and $T$, 
    such that
    \begin{equation*}
          \sum_{\iota=u,v} \|\chi^{m+1}_{\iota}\|  \le \widehat{C}_2(\tau^2+h^4)
           +   C_1 \tau \sum_{n=0}^{m} \sum_{\iota=\phi,\psi}  \|\chi_\iota^{n+1/2} \|,  \quad 0\le m \le N-1.
    \end{equation*}
\end{lemma}
\begin{proof}
Taking inner product of \eqref{err_eq:M1}  with $\chi^{n+1}_u+\chi^n_u $ and taking imaginary part of the resultant equation, we obtain from Lemma \ref{lem:operator} that
$$
  \|\chi^{n+1}_u \|-\|\chi^n_u \|\le \tau (\| \mrm D_{u,1}^{n+1}\|+\| \mrm D_{u,3}^{n+1}\|) .  
$$
  In the view of \eqref{ieq:bound_UmVm}, Lemmas \ref{lem:norm}  and  \ref{lem:trunc_errors}, we have
  \begin{equation}\label{eu}
      \|\chi^{n+1}_u \|-\|\chi^n_u \|
    \le  \f{C_1\tau}{2} \big(\|\chi_\phi^{n+1/2} \| +\|\chi_\psi^{n+1/2}\| \big)
         +\widehat{C}_0\tau(\tau^2+h^4),
  \end{equation}
and upon summing $n$ from $0$ to $m$ $(0\le m\le N-1)$,  we further obtain
\begin{equation}\label{err_eq:eu}
    \|\chi^{m+1}_u \|\le \f{C_1\tau}{2} \sum_{n=0}^{m}\big(\|\chi_\phi^{n+1/2} \| +\|\chi_\psi^{n+1/2}\| \big)
       +\widehat{C}_2 (\tau^2+h^4).
\end{equation}
 Analogously, it follows that
 \begin{equation}\label{err_eq:ev}
    \|\chi^{m+1}_v\|\le \f{C_1\tau}{2}\sum_{n=0}^{m}\big(\|\chi_\phi^{n+1/2} \| +\|\chi_\psi^{n+1/2}\| \big)
      +\widehat{C}_2(\tau^2+h^4).
\end{equation}
Thus, adding \eqref{err_eq:eu} and \eqref{err_eq:ev} together, we prove the conclusion.
\end{proof}
In particular, for the first two time levels, we have the following conclusions.
\begin{lemma}\label{lem:M_err2}
Assume that the condition of Lemma \ref{lem:M_err1} holds. There exists a positive constant $\widehat{C}_3$ related to $\widehat{C}_0$, $\widehat{C}_1$ and $C_1$, 
such that	
   \begin{equation*}
      \sum_{j=1}^2\Big( \sum_{\iota=u,v} \big(\|\chi^{j}_\iota \| 
        +\tau \| \delta_{t}\chi_\iota^j\| \big)
        + \tau \sum_{\iota=\phi,\psi}\|\chi^{j-1/2}_\iota\| \Big)
        \le \widehat{C}_3\tau(\tau^2+h^4).
    \end{equation*}
\end{lemma}
\begin{proof}
From \eqref{eu} and Lemma \ref{lem:err_ana_pre}, we see
 \begin{equation*}
	\|\chi^{1}_u \|\le \widehat{C}_0\tau(\tau^2+h^4)+C_1\tau\big(\|\chi_\phi^{1/2} \| +\|\chi_\psi^{1/2}\|\big)\le \widehat{C}_3\tau(\tau^2+h^4),
\end{equation*}
and similarly, $\|\chi^{1}_v \|\le  \widehat{C}_3\tau(\tau^2+h^4)$.
Next, taking norm $\|\cdot\|$ on the both sides of \eqref{err_eq:M3} for $n=1$, we have
 \begin{equation*}\label{lem42:e2}
  \|\chi^{3/2}_\phi \|\le \|e^{1/2}_\phi \|+\| \mathrm{r}_\phi^1\|+2\|~|u^1 |^2-|U^{1}_{\mfm} |^2 ~ \|
    \le \widehat{C}_3 (\tau^2+h^4),
\end{equation*}
where \eqref{ieq:bound_UmVm} and Lemmas  \ref{lem:trunc_errors}, \ref{lem:err_ana_pre} have been applied. Analogously, we can show 
$$
\|\chi^{3/2}_\psi \|\le  \widehat{C}_3 (\tau^2+h^4).
$$

Moreover, we see from \eqref{eu} that
\begin{equation*}
    \|\chi^{2}_u \|\le \|\chi^{1}_u \|+\widehat{C}_0\tau(\tau^2+h^4)
        +C_1\tau\sum_{n=0}^{1}\big(\|\chi_\phi^{n+1/2} \| +\|\chi_\psi^{n+1/2} \|\big)
         \le \widehat{C}_3 \tau(\tau^2+h^4).
\end{equation*}
Furthermore, we have
 \begin{equation*}
  \|\delta_t\chi^{1}_u \|=\tau^{-1} \|\chi^{1}_u \| \le \widehat{C}_3(\tau^2+h^4),~~
          \|\delta_t\chi^{2}_u \|\le \tau^{-1} (\|\chi^{1}_u \|+\|\chi^{2}_u\|)\le \widehat{C}_3 (\tau^2+h^4).
\end{equation*}
The estimates $\|\chi_v^2\|$ and $\|\delta_t\chi_v^j\|$, $j=1,2$, can be obtained in a similar way.
\end{proof}

In the following, we investigate the errors $\|\chi_\phi^{n+1/2}\|$ and $\|\chi_\psi^{n+1/2}\|$ that appeared in Lemma \ref{lem:M_err1}.
	\begin{lemma}\label{lem:M_err3}
		Assume that the condition of Lemma \ref{lem:M_err1} holds. There exists a positive constant $\widehat{C}_4$ relying on $T$, $\widehat{C}_0$ and $\widehat{C}_3$, such that	
		\begin{equation*}
			\begin{aligned}
				\sum_{\iota=\phi,\psi} \|\chi^{m+1/2}_{\iota}\| 
				\le \widehat{C}_4(\tau^2+h^4)+C_2\tau\sum_{n=2}^{m}\sum_{\iota=u,v}\big(\|\chi^{n}_\iota\|+\|\delta_{t}\chi^{n}_\iota\|\big),\quad 2\le m \le N-1.
			\end{aligned}
		\end{equation*}
	\end{lemma}
	\begin{proof}~
		Taking difference for \eqref{err_eq:M3} between $t=t_{n} $ and $t=t_{n-1}$ $(2\le n\le N-1)$, we conclude 
		\begin{equation}
			\chi^{n+1/2}_{\phi}-\chi^{n-3/2}_{\phi}=\xi_u^{n}+\mrm{r}_{\phi}^{n-1}-\mrm{r}_{\phi}^{n},\label{err_eq:phi}
		\end{equation}
		where $\xi_u^{n}:=2(|u^{n} |^2-|U_\mfm^{n} |^2  )-2(|u^{n-1} |^2-|U_\mfm^{n-1} |^2)$. Together with Lemma \ref{lem:ieq:s2} and \eqref{ieq:bound_UmVm}, it yields 
		\begin{equation} \label{err_eq:phi:e2}
			\begin{aligned}
				\|\xi^{n}_u \| 
				&\le  4\|u^n- u^{n-1}\|_\infty \|U_\mfm^{n}-u^{n} \|\\
				& \qquad +2\mrm{G_I}(U_\mfm^{n-1},U_\mfm^{n},u^{n-1},u^{n})\|U_\mfm^n-U_\mfm^{n-1}-u^n+u^{n-1} \|\\
				&\le \f{C_2\tau}{2}(\|\chi^{n}_u \|+\|\delta_{t}\chi^{n}_u\|),
			\end{aligned} 
		\end{equation}
		where
		$\|u^n- u^{n-1}\|_\infty\le \tau \|u_t\|_{L^\infty(\Omega\times I)} \le \mK\tau$ and $$\mrm{G_I} = \|U_\mfm^{n-1}\|_\infty+ \|U_\mfm^{n}\|_\infty+\|u^{n-1}- u^n\|_\infty\le 4\mK+2$$
		have been applied.     
		
		Now, taking the discrete inner product of \eqref{err_eq:phi} with $\chi^{n+1/2}_{\phi}+\chi^{n-3/2}_{\phi}$, it follows from \eqref{err_eq:phi:e2} and Lemma \ref{lem:trunc_errors} that
		$$   
		\begin{aligned}
			\|\chi^{n+1/2}_{\phi} \| -\|\chi^{n-3/2}_{\phi} \|
			\le \|\xi^{n}_u \|+\|\mrm{r}_{\phi}^{n-1}-\mrm{r}_{\phi}^{n} \| 
			\le \f{C_2\tau}{2} (\|\chi^{n}_u \|+\|\delta_{t}\chi^{n}_u\|)+\widehat{C}_0\tau^3. 
		\end{aligned} 
		$$
Therefore, summing over $n$ from $2$ to $m$ $(2\le m\le N-1)$,  and utilizing Lemma \ref{lem:M_err2}, we have
		\begin{equation*}
			\begin{aligned}
				\|\chi^{m+1/2}_{\phi} \|+\|\chi^{m-1/2}_{\phi} \|
				&\le \|\chi^{3/2}_{\phi} \| +\|\chi^{1/2}_{\phi} \| + \f{C_2\tau}{2}\sum_{n=2}^{m}(\|\chi^{n}_u \|+\|\delta_{t}\chi^{n}_u\|) +\widehat{C}_0 T\tau^2\\
				& \le \widehat{C}_4(\tau^2+h^4) +\f{C_2\tau}{2}\sum_{n=2}^{m}(\|\chi^{n}_u \|+\|\delta_{t}\chi^{n}_u\|).
			\end{aligned} 
		\end{equation*}
		Similarly, it follows that
		\begin{equation*}
			\|\chi^{m+1/2}_{\psi} \|+\|\chi^{m-1/2}_{\psi} \|  
			\le \widehat{C}_4(\tau^2+h^4) +\f{C_2\tau}{2}\sum_{n=2}^{m}(\|\chi^{n}_v \|+\|\delta_{t}\chi^{n}_v\|),
		\end{equation*}
		which implies the  conclusion.
	\end{proof}

Below, we provide the estimates of the temporal divided differences $\|\delta_t\chi^{n+1}_u\|$ and $\|\delta_t\chi^{n+1}_v\|$, which play an important role in subsequent analysis.
\begin{lemma}\label{lem:M_err4}
  Assume that the condition of Lemma \ref{lem:M_err1} holds.There are some positive constants  $\widehat{C}_5$ related to $T$, $C_1$,  $\widehat{C}_0$ and $\widehat{C}_3$  such that
\begin{equation*}
\begin{aligned}
      \sum_{\iota=u,v} \| \delta_{t}\chi^{m+1}_\iota\| 
    & \le \widehat{C}_5(\tau^2+h^4)+C_3^*\tau \sum_{n=2}^m \sum_{\iota=u,v}(\|\chi^{n+1}_\iota\|+\|\delta_{t}\chi_\iota^{n+1} \|)\\
    &\quad    +2C_3\tau \sum_{n=2}^m\sum_{\iota=\phi,\psi}\|\chi^{n-3/2}_\iota\|,	\quad 2\le m\le N-1.
\end{aligned}	
\end{equation*}
\end{lemma}
\begin{proof}
   Taking difference for \eqref{err_eq:M1} between $t=t_{m+1/2}$ and $t=t_{m-3/2}$ $(2\le m\le N-1)$, and reformulate it as 
    \begin{equation}\label{err_eq:delta2}
        \delta_{t}\chi^{m+1}_u=( \mal C_h-\mal I_h)\delta_{t}\chi^{m}_u+\mal C_h \delta_{t}\chi^{m-1}_u
                    +  \sum_{j=1}^{3}\mT_h^{-1}\mrm G_{u,j}^{m},    
    \end{equation}
where $\mrm G_{u,j}^{m}:=\mrm i\mrm D_{u,j}^{m+1}-\mrm i\mrm D_{u,j}^{m-1}$.  An application of Lemma \ref{lem:space_operator} to \eqref{err_eq:delta2} yields
\begin{equation}\label{err_eq:delta3}
    \begin{aligned}
         \delta_{t}\chi^{m+1}_u &=\f{1}{2} \big[(-1)^{m+1}\mal{I}_h+\mal{C}_h^{m}\big]\mT_h\delta_{t}\chi_u^2
                 +\f{1}{2}\big[(-1)^{m}\mal{C}_h+\mal{C}_h^{m}\big]\mT_h\delta_{t}\chi_u^1\\
                 &\quad+\f{1}{2}\sum_{n=2}^m \big[(-1)^{m-n}\mal{I}_h+\mal{C}_h^{m+1-n}\big]\sum_{j=1}^{3}\mrm G_{u,j}^{n}.
    \end{aligned} 
\end{equation}
Then, taking the discrete norm $\|\cdot\|$ on both sides of \eqref{err_eq:delta3}, and Lemma \ref{lem:OperatorNorm} implies that
 \begin{equation}\label{errdeltat}
    \| \delta_{t}\chi^{m+1}_u\|  \le 
    \|\mT_h\delta_{t}\chi_u^2 \| +\| \mT_h\delta_{t}\chi_u^1\|   +\sum_{n=2}^m\sum_{j=1}^{3}\| \mrm G_{u,j}^{n}\|. 
  \end{equation}
 
Next, we proceed to estimate the right-hand side terms of \eqref{errdeltat}. Firstly, according to \eqref{ieq:bound_UmVm} and Lemma \ref{lem:M_err2}, it follows that
\begin{equation}\label{err_eq:D1}
    \| \mrm D_{u,1}^{j}\|\le \widehat{C}_5(\tau^2+h^4), ~  \| \mrm D_{u,2}^{j}\|\le \widehat{C}_5(\tau^2+h^4),~ \| \mrm D_{u,3}^{j}\|\le \widehat{C}_5(\tau^2+h^4),~ j=1,2.
\end{equation}
Therefore, taking the discrete norm $\| \cdot\|$ on both sides of \eqref{err_eq:M1} for $n=0$ and $n=1$ respectively, we arrive at
\begin{equation}\label{errdeltat:e1}
  \begin{aligned}
   \| \mT_h\delta_{t}\chi_u^1\|  &=\Big\|\sum_{j=1}^{3} \mrm D_{u,j}^{1}\Big\|\le \widehat{C}_5(\tau^2+h^4),\\ 
   \| \mT_h\delta_{t}\chi_u^2\|  &=\Big\| -2(\mT_h-\mal I_h) \delta_{t}\chi_u^1
            +\sum_{j=1}^{3} \mrm D_{u,j}^{2}\Big\|\le \widehat{C}_5(\tau^2+h^4),
    \end{aligned}        
\end{equation}
where Lemma \ref{lem:M_err2} and  \eqref{err_eq:D1} have been applied.

Secondly, in view of the triangle inequality, we decompose $\|\mrm G_{u,1}^{n} \|$ into the following two parts 
   $$   \begin{aligned}
		\|\mrm G_{u,1}^{n} \|
		&=\|( \mfm_\phi\chi_\phi^{n+1/2}+\beta \mfm_\psi\chi_\psi^{n+1/2})\widetilde{U}_{\mfm}^{n+1/2}-( \mfm_\phi\chi_\phi^{n-3/2}+\beta \mfm_\psi\chi_\psi^{n-3/2})\widetilde{U}_{\mfm}^{n-3/2}\|\\
		&\le \|(\mfm_\phi\chi_\phi^{n+1/2}+\beta \mfm_\psi\chi_\psi^{n+1/2}-\mfm_\phi\chi_\phi^{n-3/2}-\beta \mfm_\psi\chi_\psi^{n-3/2}) \widetilde{U}_{\mfm}^{n+1/2}\|\\
		&\quad +\|( \mfm_\phi\chi_\phi^{n-3/2}+\beta \mfm_\psi\chi_\psi^{n-3/2}) (\widetilde{U}_{\mfm}^{n+1/2}-\widetilde{U}_{\mfm}^{n-3/2} )  \| =:\Theta_1+\Theta_2. 
    \end{aligned} $$
For the first part, it follows from \eqref{ieq:bound_UmVm}, \eqref{err_eq:phi}--\eqref{err_eq:phi:e2} and Lemma \ref{lem:trunc_errors} that
\begin{equation}\label{err_eq:theta1}
    \begin{aligned}
		\Theta_1
		&\le \| \widetilde{U}_{\mfm}^{n+1/2}\|_\infty\big(\|\chi_\phi^{n+1/2}-\chi_\phi^{n-3/2} \|+|\beta|\|\chi_\psi^{n+1/2}-\chi_\psi^{n-3/2}\| \big) \\
		&\le (1+|\beta|)\| \widetilde{U}_{\mfm}^{n+1/2}\|_\infty\big(\|\xi^{n}_u+\mrm{r}_{\phi}^{n-1}-\mrm{r}_{\phi}^{n} \|+\|\xi^{n}_v+\mrm{r}_{\psi}^{n-1}-\mrm{r}_{\psi}^{n} \|\big) \\
        &\le C_3\tau \sum_{\iota=u,v}\big(\|\chi^{n}_\iota \|+\|\delta_{t}\chi^{n}_\iota\|\big)+\widehat{C}_0\tau(\tau^2+h^4).
	\end{aligned} 
\end{equation}
For the second part, due to $\mfm_\iota\le 1$ $(\iota=u,v,\psi,\phi)$  and \eqref{ieq:bound_phi}, we have
\begin{equation}\label{err_eq:theta2}
    \begin{aligned}
	\Theta_2
        &=\|( \mfm_\phi\chi_\phi^{n-3/2}+\beta \mfm_\psi\chi_\psi^{n-3/2}) (\widetilde{u}^{n+1/2}-\widetilde{u}^{n-3/2}+\mfm_u\widetilde{\chi}_u^{n+1/2}-\mfm_u\widetilde{\chi}_u^{n-3/2} )\|\\
	  &\le \|\widetilde{u}^{n+1/2}-\widetilde{u}^{n-3/2} \|_\infty(\|\chi_\phi^{n-3/2}\|+|\beta |\|\chi_\psi^{n-3/2}\|) \\
      & \qquad + \f\tau2\|  \mfm_\phi\chi_\phi^{n-3/2}+\beta \mfm_\psi\chi_\psi^{n-3/2}\|_\infty \|\delta_{t}\chi_u^{n+1} +2\delta_{t}\chi_u^{n} +\delta_{t}\chi_u^{n-1} \|\\
      &\le  \f{C_1\tau}{2} ( \| \chi_\phi^{n-3/2}\| +\|\chi_\psi^{n-3/2} \|) 
       +C_1^* \tau (\|\delta_{t}\chi_u^{n+1}\|+2\|\delta_{t}\chi_u^{n} \|+\|\delta_{t}\chi_u^{n-1} \| ).
    \end{aligned} 
\end{equation}
Therefore, we conclude from \eqref{err_eq:theta1}--\eqref{err_eq:theta2} that
\begin{equation}\label{errdeltat:e3}    
   \begin{aligned}
     \|\mrm G_{u,1}^{n} \|
            &\le \widehat{C}_0\tau(\tau^2+h^4)+C_3\tau \Big(\sum_{\iota=u,v}\|\chi^{n}_\iota\|+\sum_{\iota=\phi,\psi}\|\chi^{n-3/2}_\iota\|\Big)\\
            &\quad +C_2^*\tau\big(\|\delta_{t}\chi_u^{n+1}\|+\|\delta_{t}\chi_u^{n} \|+\|\delta_{t}\chi_u^{n-1} \|+\|\delta_{t}\chi_v^{n} \|\big).
   \end{aligned} 
\end{equation}

Thirdly, denote $a^{n+1/2}:=\phi^{n+1/2}+\beta\psi^{n+1/2}$. It follows from the triangle inequality and the fact $\mfm_u\le 1$ that
\begin{equation}\label{errdeltat:e4} 
    \begin{aligned}
	\|\mrm G_{u,2}^{n} \|
		&=\|a^{n+1/2}(\mfm_u\chi_u^{n+1} +\mfm_u\chi_u^{n})- a^{n-3/2}(\mfm_u\chi_u^{n-1}+\mfm_u\chi_u^{n-2})\|\\
		&\le \|(a^{n+1/2}-a^{n-3/2} ) (\chi_u^{n+1}+\chi_u^{n}) \|+\tau\|a^{n-3/2}(\delta_{t}\chi_u^{n+1} +2\delta_{t}\chi_u^{n} +\delta_{t}\chi_u^{n-1}  ) \| \\
            &\le \|a^{n+1/2}-a^{n-3/2}\|_\infty \|\chi_u^{n+1}+\chi_u^{n}\|+\tau\|a^{n-3/2}\|_\infty \|\delta_{t}\chi_u^{n+1} +2\delta_{t}\chi_u^{n} +\delta_{t}\chi_u^{n-1}\| \\
		&\le C_1\tau(\| \chi_u^{n+1}\|+\| \chi_u^{n}\| + \|\delta_{t}\chi_u^{n+1} \|+\|\delta_{t}\chi_u^{n} \|+\|\delta_{t}\chi_u^{n-1} \|).
	\end{aligned} 
\end{equation}
   
Finally,  in light of  Lemmas \ref{lem:norm} and \ref{lem:trunc_errors}, it is easy to see that
 \begin{equation}\label{errdeltat:e5} 
     \|\mrm G_{u,3}^{n}  \|=\| \mA_h^{-1}(\mrm{r}_{u}^{n+1/2}-\mrm{r}_{u}^{n-3/2})\|\le \widehat{C}_0\tau(\tau^2+h^4). 
 \end{equation}
 
 Now, inserting the obtained estimates in \eqref{errdeltat:e1}, \eqref{errdeltat:e3}--\eqref{errdeltat:e5} into \eqref{errdeltat}, we obtain
    \begin{equation}\label{error analysis:deltateu}
			\| \delta_{t}\chi^{m+1}_u\| 
            \le \widehat{C}_5(\tau^2+h^4)+\f{C_3^*\tau}{2} \sum_{n=2}^m \sum_{\iota=u,v}(\|\chi^{n+1}_\iota\|+\|\delta_{t}\chi_\iota^{n+1} \|)+C_3\tau\sum_{n=2}^m \sum_{\iota=\phi,\psi}\|\chi^{n-3/2}_\iota\|. 		
    \end{equation}
Similarly, we can obtain
    \begin{equation}\label{error analysis:deltatev}
		\| \delta_{t}\chi^{m+1}_v\| \le \widehat{C}_5(\tau^2+h^4)+
           \f{C_3^*\tau}{2} \sum_{n=2}^m \sum_{\iota=u,v}(\|\chi^{n+1}_\iota\|+\|\delta_{t}\chi_\iota^{n+1} \|)
            +C_3\tau\sum_{n=2}^m \sum_{\iota=\phi,\psi}\|\chi^{n-3/2}_\iota\|. 	
    \end{equation}
    Thus, combining  \eqref{error analysis:deltateu}  and \eqref{error analysis:deltatev} we complete the proof.
\end{proof}

Next, we provide the convergence results of the auxiliary error equations \eqref{err_eq:M1}--\eqref{err_eq:M4}.

\begin{theorem}\label{thm:M_err5} Under the regularity condition \eqref{assu:Regular}, there exists a positive constant $C^*$, independent of $\tau$ and $h$ but related to $\mK^*$, such that	
\begin{equation*}\label{err_eq:chi_L2_1}
    \max_{1\le n\le N}\Big(  \|\chi_u^{n}  \| +\|\chi_v^{n} \| 
      +  \|\chi_\phi^{n-1/2}  \|
         +\|\chi_\psi^{n-1/2}\| +\|\delta_{t}\chi_u^{n} \| +\|\delta_{t}\chi_v^{n}\|\Big)
          \le C^* \big(\tau^2 + h^4\big),
\end{equation*}
    for $\tau\le \tau_1:=\f{1}{2C_3^*} $. 
\end{theorem}
\begin{proof}
    For $1\le n\le N-1$, denote
$$
    \gamma^n:=\|\chi_u^n \|+\|\chi_v^n \| +\|\chi^{n-1/2}_{\phi} \| +\|\chi^{n-1/2}_{\psi} \| +\| \delta_{t}\chi^{n}_u\| +\| \delta_{t}\chi^{n}_v\|.
$$
Then, Lemma \ref{lem:M_err2} shows that
\begin{equation}\label{err_eq:chi_L1a}
  \gamma^{1}\le C(\tau^2+h^4),~~ 	\gamma^{2}\le C(\tau^2+h^4).
\end{equation}
Moreover, combining Lemmas \ref{lem:M_err1}, \ref{lem:M_err3}--\ref{lem:M_err4} together,  we have
\begin{equation}\label{err_eq:chi_L1b}
  \gamma^{m+1}\le C_3^*\tau\gamma^{m+1}+C^*\tau\sum_{n=2}^{m}\gamma^n +C(\tau^2+h^4),~~2\le m \le N-1, 
\end{equation}
for some positive constant $C^*$. 
Thus, taking $\tau\le\tau_1$, an application of the discrete Gronwall’s inequality to \eqref{err_eq:chi_L1b} yields
\begin{equation*} 
  \gamma^{m+1}\le C^*(\tau^2+h^4),~~2\le m \le N-1,		
\end{equation*}
which, together with \eqref{err_eq:chi_L1a} further leads to \eqref{err_eq:chi_L2_1}.
\end{proof}

 We note that the constant $C^*$ in Theorem \ref{thm:M_err5} is related to the to-be-determined constant $\mK^*$. In addition, to obtain the unconditional optimal-order error estimates of the original LRCD scheme, we still need to fix the uniform bounds for $\chi_u^n$, $\chi_v^n$, $\chi^{n-1/2}_{\phi}$ and $\chi^{n-1/2}_{\psi}$ to ensure $\mfm_u=\mfm_v =\mfm_\phi=\mfm_\psi \equiv 1$ in \eqref{cut_off:e1}.  The results are
presented in the following theorem.
 \begin{theorem}\label{thm:M_err6} Let 
$\{\chi_u^{n+1}, \chi_v^{n+1},\chi_\phi^{n+1/2},\chi_\psi^{n+1/2}\}\in \mbb{C}_h^p\times \mbb{C}_h^p\times \mathbb{R}_h^p\times \mathbb{R}_h^p$ be solutions to the auxiliary error equations \eqref{err_eq:M1}--\eqref{err_eq:M4}. Under the regularity condition \eqref{assu:Regular}, if $\tau$ and $h$ are small enough and $\tau\le \tau_1 $, there exists a positive constant $C_\ddagger$, independent of $\tau$, $h$ and $\mK^*$,  only related to $\mK$, $T$  and the exact solutions in corresponding norms, such that	
    \begin{equation*}
        \begin{aligned}
			&\max_{1\le n\le N}\|\chi_u^n\|_\infty\le 1-\epsilon, ~~\max_{1\le n\le N}\|\chi_\phi^{n-1/2}\|_\infty\le \mK^*-\epsilon,\\
			&\max_{1\le n\le N}\|\chi_v^n\|_\infty\le 1-\epsilon, ~~\max_{1\le n\le N}\|\chi_\psi^{n-1/2}\|_\infty\le \mK^*-\epsilon,
	\end{aligned}
    \end{equation*}
  for $\epsilon$ given in Section \ref{sec:AuEE} and $\mK^*:=C_\ddagger+2\epsilon$.
 \end{theorem}
  \begin{proof} 
To derive the uniform bounds of $ \chi_u^n $, $ \chi_v^n $, $ \chi^{n-1/2}_{\phi}$ and $\chi^{n-1/2}_{\psi}$ without any coupling mesh condition, we divide the following discussion into two separate cases regarding the relationship between $\tau$ and $h$. First, suppose that $\tau>h$. Note that \eqref{err_eq:M1} shows  
$$   \begin{aligned}
		\| \Lambda_h\widetilde{\chi}_u^{n+1/2}\| &= |\kappa^{-1}|\Big\| \mA_h \delta_{t}\chi^{n+1}_u +\sum_{j=1}^{3} \mA_h\mrm D_{u,j}^{n+1}\Big\|,\\
		&\le |\kappa^{-1}| \big(\|\delta_{t}\chi^{n+1}_u \|+\|\chi_\phi^{n+1/2} \|+\|\chi_\psi^{n+1/2} \|+\|\chi_u^{n+1} \|+ \|\chi_u^{n} \|+\|\mrm r_u^{n+1/2}\|\big) \\
		&\le C^* \big(\tau^2 + h^4\big),
\end{aligned} $$ 
where Lemma \ref{lem:norm}, \eqref{ieq:bound_UmVm} and Theorem \ref{thm:M_err5} have also been used. 
Thus, we conclude from Lemma \ref{lem:dis_embed} that
  \begin{equation}\label{err_eq:inf_uv_sob}
        \|\widetilde{\chi}_u^{n+1/2}\|_\infty \le C^{*}(\tau^2+h^4).
\end{equation}
Similarly, we have $\|\widetilde{\chi}_v^{n+1/2}\|_\infty\le C^{*}(\tau^2+h^4).$
 Back to the estimate $\Theta_2$ and together with Theorem \ref{thm:M_err5}, we have
 $$
    \begin{aligned}
	\Theta_2
	  &\le\big(\|\widetilde{u}^{n+1/2}-\widetilde{u}^{n-3/2} \|_\infty+ \|\widetilde{\chi}_u^{n+1/2}  -\widetilde{\chi}_u^{n-3/2} \|_\infty\big) \|\mfm_\phi\chi_\phi^{n-3/2}+\beta \mfm_\psi\chi_\psi^{n-3/2}\|\\
      &\le  \f{C_1\tau }{2} ( \| \chi_\phi^{n-3/2}\| +\|\chi_\psi^{n-3/2} \|) 
       +C^*(\tau^2+h^4)( \| \chi_\phi^{n-3/2}\| +\|\chi_\psi^{n-3/2} \|) \\
       &\le C_1\tau ( \| \chi_\phi^{n-3/2}\| +\|\chi_\psi^{n-3/2} \|), 
    \end{aligned} 
$$
for sufficiently small $\tau$ and $h$, which together with the estimate $\Theta_1$ in \eqref{err_eq:theta1} further leads to a $\mK^*$-independent estimate for $\mrm G_{u,1}^{n}$ that
$$  \begin{aligned}
    \|\mrm G_{u,1}^{n} \|
            &\le \widehat{C}_0\tau(\tau^2+h^4)+C_3\tau \Big(\sum_{\iota=u,v}(\|\chi^{n}_\iota\|+\|\delta_{t}\chi_\iota^{n} \|)+\sum_{\iota=\phi,\psi}\|\chi^{n-3/2}_\iota\|\Big).
   \end{aligned}
$$
Then, combined with the original estimates in \eqref{errdeltat:e4}--\eqref{errdeltat:e5}, we conclude 
$$
 \begin{aligned}
   & \|\delta_{t}\chi^{m+1}_u\|+\| \delta_{t}\chi^{m+1}_v\| \\
   &\le\widehat{C}_5(\tau^2+h^4)+ C_4 \tau \sum_{n=1}^m
     \Big(\sum_{\iota=u,v}  \big(\|\delta_{t}\chi_{\iota}^{n+1} \| 
           +\|\chi^{n+1}_{\iota} \|\big) 
     +\sum_{\iota=\phi,\psi}\|\chi^{n-3/2}_\iota\|\Big).
 \end{aligned}
$$
Therefore, for $\tau\le\tau_1$,  along with Lemmas \ref{lem:M_err1}--\ref{lem:M_err3} and  the discrete Gronwall’s inequality again, we have
\begin{equation}\label{err_eq:conv2}
    \max_{1\le n\le N} \Big( \|\chi_u^{n}  \|+\|\chi_v^{n}  \| + \|\chi_\phi^{n-1/2} \|
       +\|\chi_\psi^{n-1/2}\|+\|\delta_{t}\chi_u^{n}  \|+\|\delta_{t}\chi_v^{n}\|\Big)
      \le C_\sharp  (\tau^2+h^4),
 \end{equation}
 which is similar to the estimates in Theorem \ref{thm:M_err5}, but the constant $C_\sharp  $ is independent of $\tau$, $h$ and $\mK^*$, and only related to $\mK$ and the exact solutions in corresponding norms. 
Furthermore, similar to the derivation of \eqref{err_eq:inf_uv_sob}, we conclude from \eqref{err_eq:conv2} to obtain the following updated result 
$$
  \|\widetilde{\chi}_u^{n+1/2}\|_\infty =\f12 \|\chi_u^{n+1}+\chi_u^{n}\|_\infty\le C_\sharp (\tau^2+h^4),$$
which leads to
\begin{equation}\label{err_eq:u:e1}
\begin{aligned}
    \|\chi_u^{m+1}\|_\infty
    &\le \sum_{n=0}^{m} (\|\chi_u^{n+1}\|_\infty- \|\chi_u^{n}\|_\infty)
      \le \sum_{n=0}^{m}\|\chi_u^{n+1}+\chi_u^{n}\|_\infty\\
    &\le C_\sharp \tau^{-1}(\tau^2+h^4)
    \le C_\sharp (\tau+h^{3}) \le1-\epsilon,~0\le m\le N-1,
\end{aligned}
\end{equation}
for sufficiently small $\tau$ and $h$, and $\tau>h$.

In addition, taking the discrete norm $\|\cdot\|_\infty$ on both sides of \eqref{err_eq:M3} and summing over $n$ from 1 to $m$, it follows from  Lemma \ref{lem:trunc_errors} and \eqref{err_eq:u:e1} that
		$$
		\begin{aligned}
			\|\chi^{m+1/2}_{\phi}\|_\infty-\|\chi^{1/2}_{\phi}\|_\infty&\le\sum_{n=1}^{m}\|\chi^{n+1/2}_{\phi}+\chi^{n-1/2}_{\phi}\|_\infty\\
			&\le \sum_{n=1}^{m}\Big[2\, (\|u^{n}\|_\infty+\|U_\mfm^{n} \|_\infty)\|\chi_u^{n}\|_\infty+\|\mrm{r}_{\phi}^{n}\|_\infty\Big]\\
			&\le C_\sharp(4\mK+2) \sum_{n=0}^{m} (\tau+h^{3}) +\widehat{C}_0T\tau,
		\end{aligned}	
		$$
		which together with Lemma \ref{lem:err_ana_pre} leads to 
		\begin{equation}\label{err_eq:phi:e1}
			\|\chi^{m+1/2}_{\phi}\|_\infty\le C_\ddagger\tau^{-1}(\tau+h^{3})\le C_\ddagger(1+h^{2})\le C_\ddagger+\epsilon =\mK^*-\epsilon,
		\end{equation}
for fixed $\mK^*:=C_\ddagger+2\epsilon$, where the constant $C_\ddagger$ is mesh-independent and only related to $C_\sharp$, $\mK$, $T$, $\widehat{C}_0$ and $\widehat{C}_1$.  

Second, if $\tau\le h $, we have ${\tau^2}/{h}\le \tau$, and from Theorem \ref{thm:M_err5} and the inverse inequality, it is easy to check
\begin{equation}\label{err_eq:u_phi:e2}
\begin{aligned}
    &\|\chi_u^{n}\|_\infty\le C_{inv} h^{-1} \|\chi_u^{n}\|
    \le C_{inv}C^*\big(\tau+h^{3}\big)\le 1-\epsilon,\\
    &\|\chi_\phi^{n+1/2}\|_\infty\le C_{inv}C^*h^{-1} \big(\tau^2+h^4\big)\le C_{inv}C^* \big(\tau+h^{3}\big)\le \mK^*-\epsilon,
\end{aligned}
\end{equation}
for sufficiently small $\tau$ and $h$.

In summary, without any coupling mesh condition between $\tau$ and $h$, we conclude from \eqref{err_eq:u:e1}--\eqref{err_eq:u_phi:e2} that 
$$
  \|\chi_u^{n}\|_\infty\le 1-\epsilon, \quad \|\chi_\phi^{n+1/2}\|_\infty\le \mK^*-\epsilon.
$$
Similarly, we can show $ \|\chi_v^{n}\|_\infty\le 1-\epsilon$ and
$\|\chi^{n+1/2}_{\psi}\|_\infty\le \mK^*-\epsilon$. Thus, the proof is completed.
\end{proof}

Finally, we present the main theorem of this paper, establishing the uniform bounds of the numerical solutions and unconditional optimal-order error estimates of the LRCD scheme.
\begin{theorem}\label{thm:Convergence} Let $\{U^{n+1},V^{n+1},\Phi^{n+3/2},\Psi^{n+3/2} \}\in \mbb{C}_h^p\times \mbb{C}_h^p\times \mathbb{R}_h^p\times \mathbb{R}_h^p$ be solutions of the LRCD scheme \eqref{LRCD:e1a}--\eqref{LRCD:e3c}. Under the conditions of Theorem \ref{thm:M_err5}, the solutions of the LRCD scheme are uniformly bounded in the sense that
	\begin{equation}\label{uniform:bound}
		\begin{aligned}
			&\max_{1\le n\le N}\|U^n\|_\infty\le  \mK+1,~\max_{1\le n\le N}\|\Phi^{n-1/2}\|_\infty\le  \mK+\mK^*,\\
			&\max_{1\le n\le N}\|V^n\|_\infty\le  \mK+1,~	\max_{1\le n\le N}\|\Psi^{n-1/2}\|_\infty\le  \mK+\mK^*.
		\end{aligned}
	\end{equation}
Furthermore, there exist a	positive constant $C$, independent of $\tau$ and $h$, such that	
\begin{equation}\label{Err:bound:un}
 \max_{1\le n\le N} \big( \|U^{n}-u^{n}\|_1 + \|V^{n}-v^{n} \|_1
                +\|\Phi^{n-1/2}-\phi^{n-1/2} \|
                +\|\Psi^{n-1/2}-\psi^{n-1/2} \|\big)
            \le C(\tau^2+h^4).
\end{equation}
\end{theorem}
\begin{proof} According to Theorem \ref{thm:M_err6}, we see 
	$
	\mfm_u = \mfm_v =\mfm_\phi = \mfm_\psi \equiv 1
	$
in \eqref{cut_off:e1}, which means that the auxiliary error equations \eqref{err_eq:M1}--\eqref{err_eq:M4} are totally equivalent to the original error equations \eqref{err_eq:e1}--\eqref{err_eq:e4}. Hence, an application of Theorem  \ref{thm:M_err5}  yields
    \begin{equation}\label{ieq:L2_est}
        \max_{1\le n\le N}( \|e_u^{n}  \|+\|e_v^{n}  \|+\|e_\phi^{n-1/2}  \|+\|e_\psi^{n-1/2}  \|+\|\delta_{t}e_u^{n}  \|+\|\delta_{t}e_v^{n}\|)\le C(\tau^2+h^4),
    \end{equation}
    and moreover, Theorem  \ref{thm:M_err6} implies that
    \begin{equation}\label{ieq:bound_e}
        \begin{aligned}
			&\max_{1\le n\le N}\|e_u^n\|_\infty\le 1-\epsilon,~\max_{1\le n\le N}\|e_\phi^{n-1/2}\|_\infty\le \mK^*-\epsilon,\\
			&\max_{1\le n\le N}\|e_v^n\|_\infty\le 1-\epsilon,~\max_{1\le n\le N}\|e_\psi^{n-1/2}\|_\infty\le \mK^*-\epsilon.
         \end{aligned}
    \end{equation}
Thus, the uniform bounds \eqref{uniform:bound} are proved following \eqref{ieq:bound_e} and assumption \eqref{bound_assu}.

Furthermore, acting the operator $\mA_h$ on both sides of \eqref{err_eq:e1}, and then taking inner product of the resulting equation with $\delta_te^{n+1}_u$, and finally choosing the real part, we arrive at 
$$
  -(\Lambda_he^{n+1}_u,e^{n+1}_u) + (\Lambda_he^{n}_u,e^{n}_u)
  \le C\tau \Big(\sum_{j=1}^3\|\Gamma_{u,j}^{n+1}\|^2+\|\delta_te^{n+1}_u\|^2\Big),
 $$
where Lemma \ref{lem:norm} and the Cauchy-Schwarz inequality have been used. Then, summing over $n$ from $0$ to $m$ for $0\le m\le N-1$, and following the definitions of $\Gamma_{u,j}^{n+1}$, Lemmas \ref{lem:H1} and \ref{lem:trunc_errors}, the uniform bounds \eqref{uniform:bound} and error estimates in \eqref{ieq:L2_est}, we have
\begin{equation}\label{err_eq:H1}
    \begin{aligned}
        \f23 |e^{m+1}_u|_1^2
        &\le  -(\Lambda_he^{m+1}_u,e^{m+1}_u)\\
       & \le C\tau\sum_{n=0}^m \Big(\sum_{j=1}^3\|\Gamma_{u,j}^{n+1}\|^2+\|\delta_te^{n+1}_u\|^2\Big)\\
   &\le C(\tau^2+h^4)^2 \qquad \Longrightarrow \|e^{m+1}_u\|_1\le C(\tau^2+h^4).
   \end{aligned}
\end{equation}
Following the same lines of proof for \eqref{err_eq:H1}, we can also show that $\|e^{m+1}_v\|_1\le C(\tau^2+h^4).$ These together with \eqref{ieq:L2_est} imply the unconditional optimal-order error estimates \eqref{Err:bound:un} for the LRCD scheme.
\end{proof}

\begin{remark} We remark that all lemmas in Section \ref{sec:LRCD} also hold for the 1D and 3D models, but may have different constants for lower and upper bounds. For example, in 3D case, the lower bound in Lemma \ref{lem:H1} should be $\f13|w|_1^2$. 
Therefore, using the same arguments as in Sections \ref{sec:AuEE}--\ref{sec:err_est}, and by discussing the relationship between $\tau$ and $h^{d/2}$, and applying the inverse inequality $\|w\|_\infty\le C_{inv}h^{-d/2}\|w\|$ in \eqref{err_eq:u_phi:e2} for $\tau\le h^{d/2}$ $ (d \le 3)$, we can derive the unconditional optimal-order error estimates in Theorem \ref{thm:Convergence} for the 1D and 3D spaces. In Section \ref{sec:num}, we also provide a numerical example to verify the convergence order for the 3D model. 
\end{remark}

\section{Numerical examples}\label{sec:num}
In this section, we present several numerical examples to validate our theoretical conclusions, including the convergence test and conservation validation for long-term simulations of the LRCD scheme. 
All computations are performed using Matlab R2021b on a Dell Precision 7920 Tower workstation with the configuration: Intel(R) Xeon(R) Gold 6258R CPU @ 2.70GHz  and 256 GB RAM.

\subsection{Convergence test}
In this subsection, we demonstrate the convergence orders stated in Theorem \ref{thm:Convergence} for model problem in the rectangular domain ($d = 2$) and the cuboid domain ($d = 3$). For the tests below, we always set $N = M_x^2 = M_y^2 (= M_z^2) =: M^2$ to observe the fourth-order convergence with respect to  the spatial grid size. Consider the following $d$-dimensional CNLS system
$$
 \left\{
 \begin{aligned}
     \mrm i u_{t}+\Delta u+(|u|^{2}+|v|^{2}) u =f_1, ~~ \text{in} ~(0,2\pi)^d\times(0, 1],\\
     \mrm i v_{t}+\Delta v+(|v|^{2}+|u|^{2}) v =f_2,  ~~ \text{in} ~(0,2\pi)^d\times(0, 1],
      \end{aligned}
 \right.
$$
with given source functions $f_1$ and $f_2$ such that the manufactured exact solutions are
\begin{itemize}
    \item $d=2$: ~~$u(x,y,t)=e^{\mrm it}\cos x\sin y, ~~v(x,y,t)=0.5e^{\mrm it}\sin x\sin y;
$
    \item $d=3$: ~~$u(x,y,z,t)=e^{\mrm it}\cos x \sin y \cos z, ~~v(x,y,z,t)=0.5e^{\mrm it}\sin x \sin y \cos z. 
$
\end{itemize}

The errors and convergence orders (C.O.) are listed in Tables \ref{tab:dis_uv2D}--\ref{tab:dis_phipsi3D}. As can be observed, whether for the two or three space dimensions, the LRCD scheme can achieve the optimal convergence order of $\mO(\tau^2 + h^4)$, specifically, for the primal variables $U$ and $V $ in the $H^1$-norm at the temporal nodes, and for the relaxation variables $\Phi$ and $\Psi$ in the $L^2$-norm at the intermediate nodes. These findings align perfectly with the theoretical predictions of Theorem \ref{thm:Convergence}.
\begin{table}[h]
    \centering
    \caption{Errors and convergence orders of $U$ and $V$ at $T=1$ for 2D model}
    \begin{tabular}{ccccccccc}
	\toprule
    $M$ & $\|U^{N}-u^{N}\|$ & C.O. & $\|U^{n}-u^{n}\|_1$ & C.O.   & $\|V^{N}-v^{N}\|$ & C.O. & $\|V^{N}-v^{N}\|_1$ & C.O. \\
		\midrule
          8	&1.1791E-02	&{--} 	&1.9367E-02	&{--}   &4.6414E-03	&{--} &7.3499E-03	&{--}	\\  
          16	&7.2101E-04	&4.03	&1.2297E-03	&3.98   &2.8500E-04	&4.03 &4.7265E-04	&3.96\\
          32	&4.4855E-05	&4.01	&7.7482E-05	&3.99   &1.7732E-05	&4.01 &3.0141E-05	&3.97\\
          64	&2.8003E-06	&4.00 	&4.8606E-06	&3.99   &1.1070E-06	&4.00 &1.9042E-06	&3.98 \\
		\bottomrule
	\end{tabular}
    \label{tab:dis_uv2D}
\end{table}
\begin{table}[H]
    \centering
    \caption{Errors and convergence orders of $\Phi$ and $\Psi$  at $t=T-\tau/2$ for 2D model}
\begin{tabular}{ccccc}
    \toprule
      $M$ & $\|\Phi^{N-1/2}-\phi^{N-1/2}\|$ & C.O.  & $\|\Psi^{N-1/2}-\psi^{N-1/2}\|$ & C.O. \\
    \midrule
        8	&1.5688E-02	&{--}	    &3.2626E-03	&{--}\\
        16	&9.4294E-04	&4.06	&1.9866E-04	&4.04\\
        32	&5.8520E-05	&4.01 	&1.2342E-05	&4.01\\
        64	&3.6512E-06	&4.00 	&7.7023E-07	&4.00 \\
    \bottomrule
    \end{tabular}
    \label{tab:dis_phipsi2D}
\end{table}
\begin{table}[h]
    \centering
    \caption{Errors and convergence orders of $U$ and $V$ at $T=1$ for 3D model}
      \begin{tabular}{ccccccccc}
		\toprule
    $M$ & $\|U^{N}-u^{N}\|$ & C.O. & $\|U^{n}-u^{n}\|_1$ & C.O.   & $\|V^{N}-v^{N}\|$ & C.O. & $\|V^{N}-v^{N}\|_1$ & C.O. \\
		\midrule
		8	&2.1567E-02	&{--}	    &4.0859E-02	&{--}   &8.4216E-03	&{--}		&1.5548E-02	&{--}\\
            12	&4.1805E-03	&4.05	&8.1569E-03	&3.97	&1.6389E-03	&4.04	&3.1203E-03	&3.96\\
            18	&8.2031E-04	&4.02	&1.6235E-03	&3.98	&3.2173E-04	&4.02	&6.2443E-04	&3.97\\
            27	&1.6159E-04	&4.01	&3.2205E-04	&3.99	&6.3381E-05	&4.01	&1.2450E-04	&3.98\\
		\bottomrule
	\end{tabular}
    \label{tab:dis_uv3D}
\end{table}
\begin{table}[H]
    \centering
    \caption{Errors and convergence orders of $\Phi$ and $\Psi$  at $t=T-\tau/2$ for 3D model}
\begin{tabular}{ccccc}
    \toprule
      $M$ & $\|\Phi^{N-1/2}-\phi^{N-1/2}\|$ & C.O.  & $\|\Psi^{N-1/2}-\psi^{N-1/2}\|$ & C.O. \\
    \midrule
       8	&2.9480E-02	   &{--}	   &5.4225E-03 &{--}\\
      12	 &5.6701E-03	&4.07	&1.0544E-03	&4.04\\
      18	 &1.1114E-03	&4.02	&2.0713E-04	&4.01\\
      27	 &2.1886E-04	&4.01	&4.0822E-05	&4.01\\
    \bottomrule
    \end{tabular}
    \label{tab:dis_phipsi3D}
\end{table}

\subsection{Conservation validation}
In this subsection, we conduct a numerical example to verify the long-term conservation properties established in Theorem \ref{thm:conservation}. Let us choose the parameter values $\beta=1.5$, $\kappa=0.5$ and the initial conditions
$$
u_0(x,y)=0.5e^{-x^2-y^2}, ~~v_0(x,y)=0.3e^{-(x-5)^2-(y-5)^2},\quad \text{in}~ \Omega=(-10,10)^2,
$$
for model \eqref{model:e1a}--\eqref{model:e1b}.

\begin{figure} [H]
    \centering  
   \vspace{-0.3cm} 
       \includegraphics[width=0.45\linewidth]{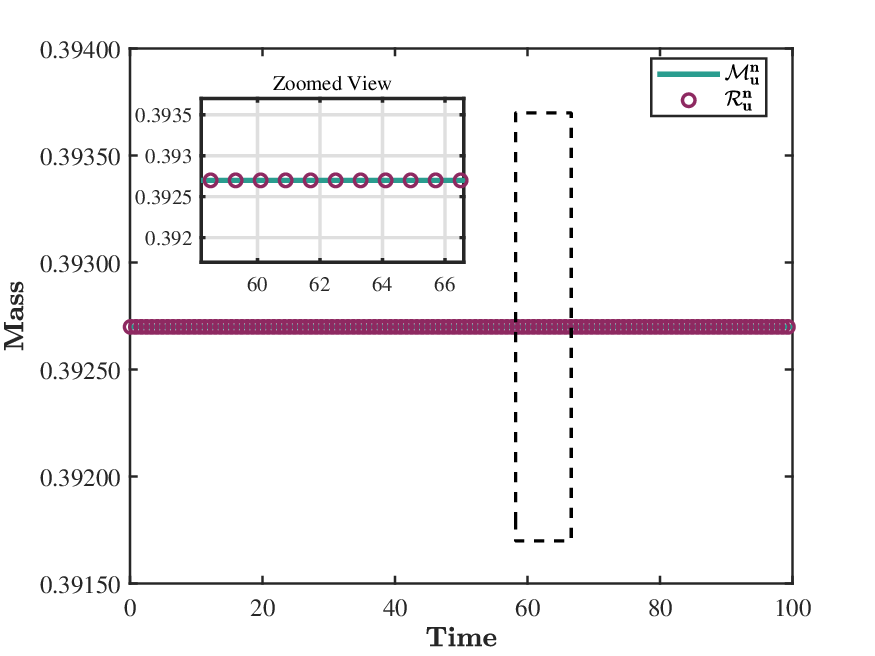}
       \includegraphics[width=0.45\linewidth]{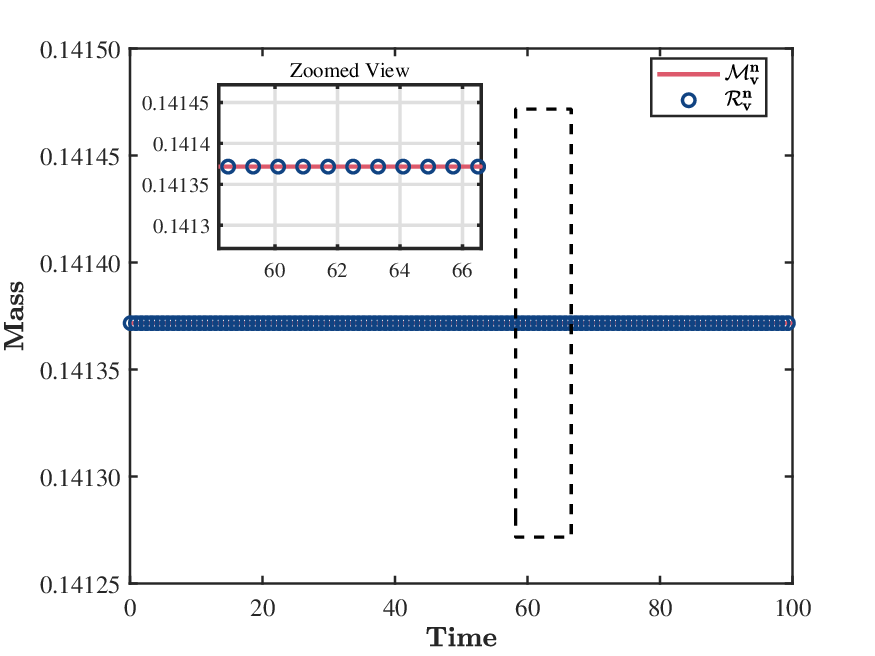}
    \caption{Evolution of mass with respect to time.}    \label{fig:mass}
\end{figure}

\begin{figure}[H]
    \centering  
    \vspace{-0.3cm} 
       \includegraphics[width=0.45\linewidth]{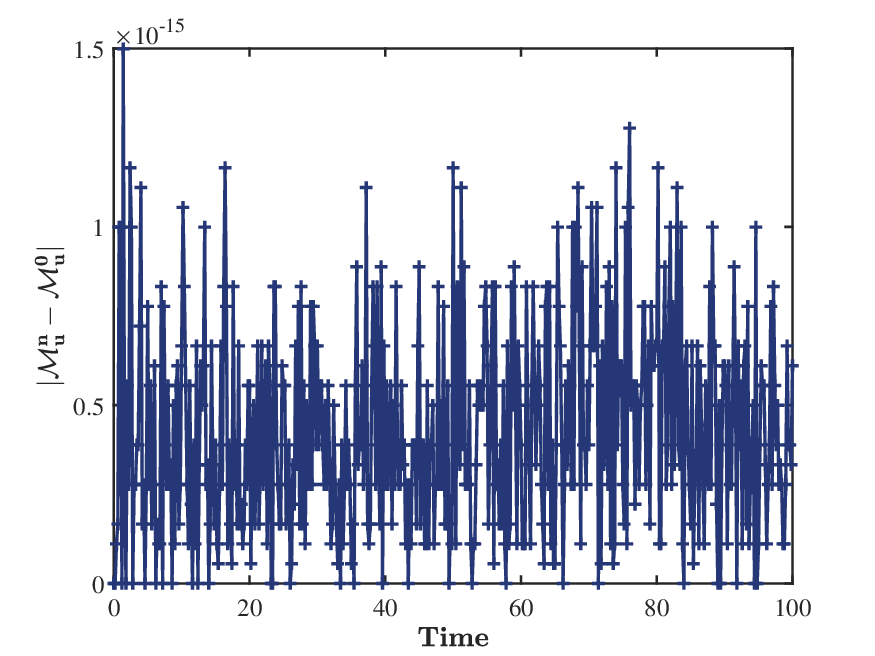}
       \includegraphics[width=0.45\linewidth]{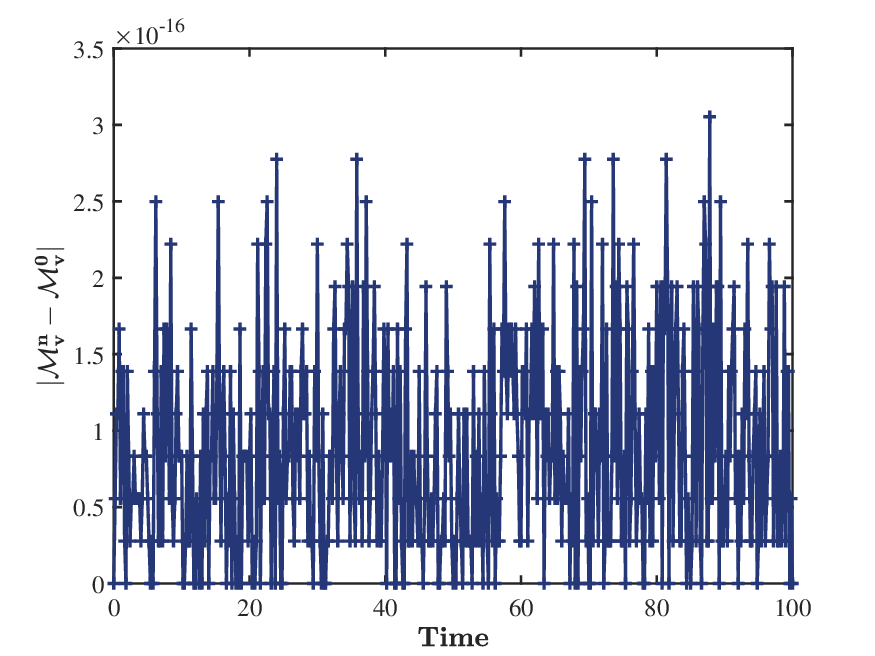}
    \caption{Evolutions of mass errors $|\mal M_u^n-\mal M_u^0|$ (left) and  $|\mal M_v^n-\mal M_v^0|$ (right) with respect to time. }
    \label{fig:err_M_u_M_v}
\end{figure}

\begin{figure}[H]
    \centering  
    \vspace{-0.3cm} 
	 \includegraphics[width=0.45\linewidth]{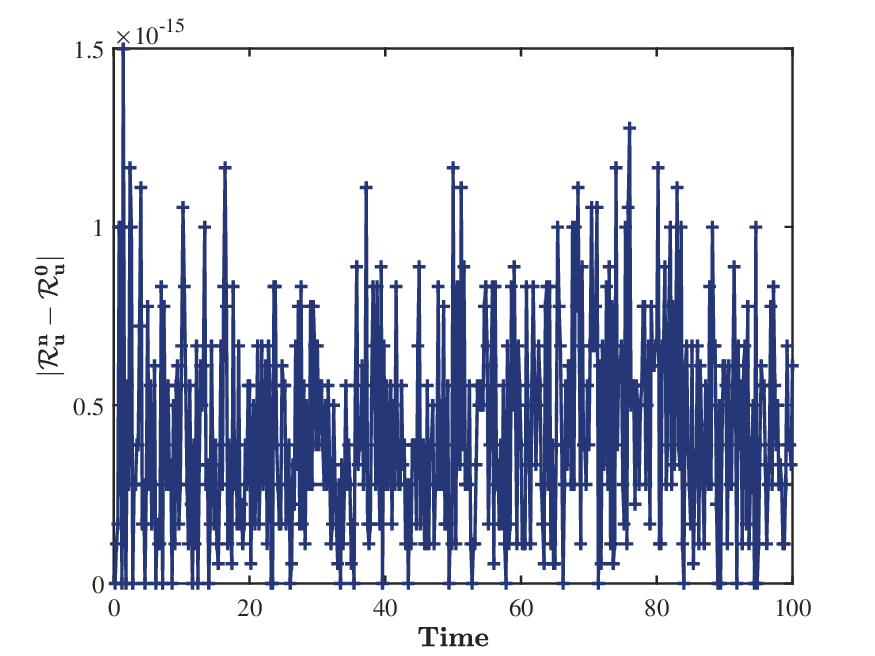}
	\includegraphics[width=0.45\linewidth]{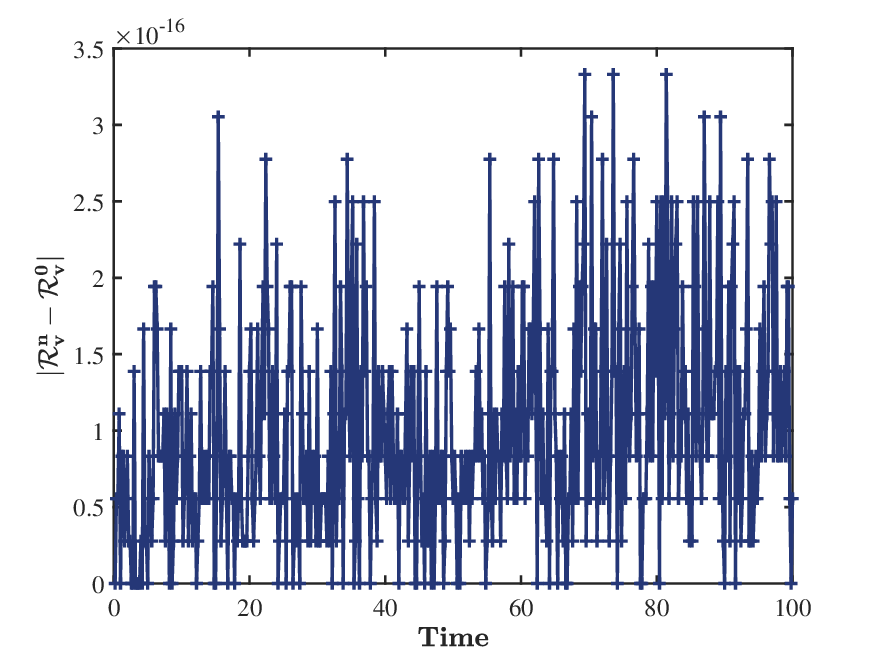}
    \caption{Evolution of  mass errors  $|\mal R_u^n-\mal R_u^0|$  (left) and $|\mal R_v^n-\mal R_v^0|$ (right) with respect to time. }    \label{fig:err_R_u_R_v}
\end{figure}

\begin{figure}[H]
    \centering  
    \vspace{-0.3cm} 
		\includegraphics[width=0.45\linewidth]{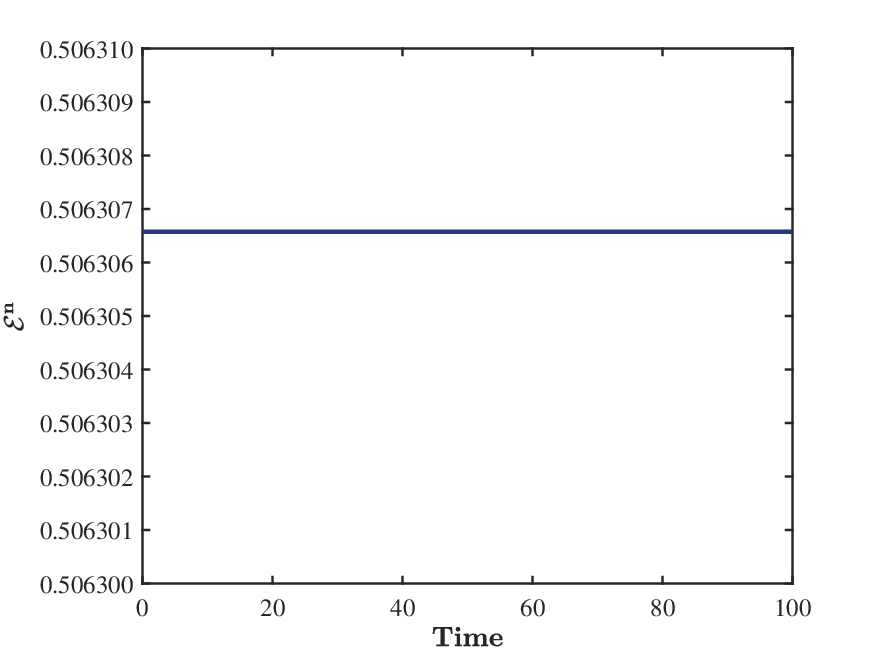}
		\includegraphics[width=0.45\linewidth]{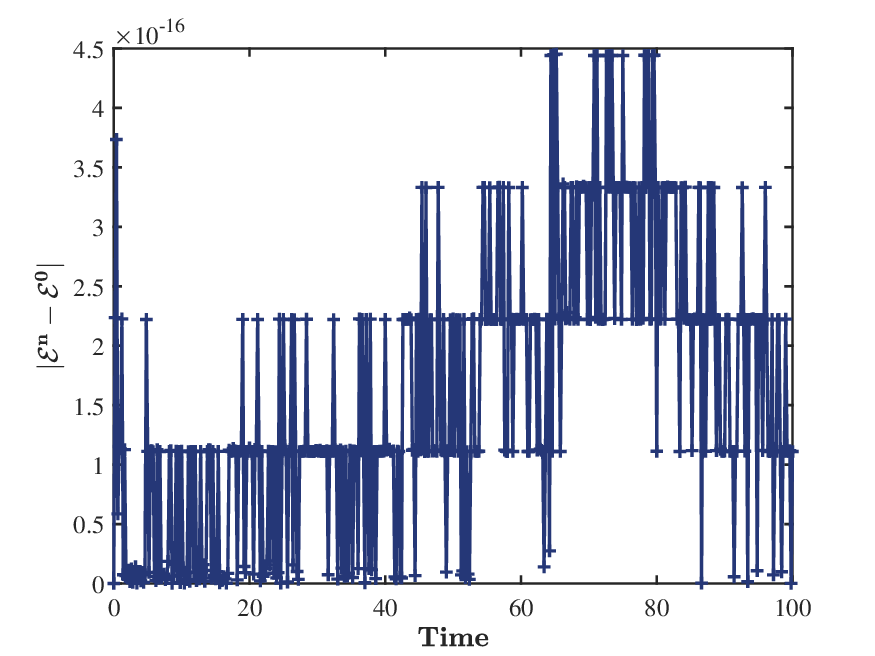}
    \caption{Evolutions of energy (left) and absolute error (right) with respect to time.}
    \label{fig:energy}
\end{figure}

The evolutions of mass, energy and their absolute errors until $T=100$ are shown in Figures \ref{fig:mass}--\ref{fig:energy} with fixed $\tau=h=0.2$.  Figure \ref{fig:mass} indicates that the total mass remains constant over time and the two different forms of mass are almost equal, i.e., $\mM_u^{n}=\mal R_u^{n}$, $\mM_v^{n}=\mal R_v^{n}$.  Moreover, as revealed in Figures \ref{fig:err_M_u_M_v}--\ref{fig:err_R_u_R_v}, the absolute errors for $\mal M_u^n$, $\mal M_v^n$, $\mal R_u^n$ and $\mal R_v^n$ are remarkably small, at least in the order of magnitude of $10^{-15}$. In addition, Figure \ref{fig:energy} confirms that the total energy is also well preserved and the absolute error of $\mE^n$ can reach the order of magnitude of $10 ^ {-16} $. It is worth noting that the absolute errors shown in Figures \ref{fig:err_M_u_M_v}--\ref{fig:energy} exhibit no significant growth over time, underscoring the stability of the LRCD scheme in long-term simulations.

\subsection{Collision of two solitons}
This subsection investigates the effectiveness and reliability of the developed LRCD scheme in simulating the dynamics of the interaction between two solitons in one space dimension \cite{WTC2008,SWW2017, WTC2010}. Consider the following one-dimensional CNLS system
$$
 \begin{cases}
     \mrm i u_{t}+ u_{xx}+(|u|^{2}+\beta|v|^{2}) u =0,~&(x,t)\in (-40,40)\times[0,80], \\
     \mrm i v_{t}+ v_{xx}+(|v|^{2}+\beta|u|^{2}) v =0,~&(x,t)\in (-40,40)\times[0,80],
 \end{cases}
$$
with the initial conditions
$$
u_{0}(x)=\sqrt{2}r_{1} \sech (r_{1}x+0.5\lambda_{1})e^{\mrm i\alpha_{1}x},
~v_{0}(x)=\sqrt{2}r_{2}\sech(r_{2}x+0.5\lambda_{2})e^{\mrm i\alpha_{2}x}.
$$

In the following simulation, we consider the parameter configuration with $r_1=r_2=1$,  $\lambda_1=-\lambda_2=18$ and $\alpha_1=-\alpha_2=\f{\alpha}{4}$. As discussed in Refs. \cite{WTC2008,SWW2017}, soliton interactions can be classified into two distinct regimes: elastic and inelastic collisions. In an elastic collision, the solitons recover their original waveforms after interaction, which occurs when $\beta=0$ or $\beta=1$. Conversely, the inelastic collision leads to permanent waveform modifications, involving reflection and mutual entanglement. Setting $\tau=0.01$ and $h=0.1$, we simulate three different collision scenarios with different parameters as follows:
\begin{itemize}
       \item Elastic collision: 
             $\alpha=1,~\beta=1;$
       \item Inelastic collision (reflection):
             $\alpha=1.15,~\beta=2/3;$
        \item Inelastic collision (mutual entanglement):             $\alpha=1.05,~\beta=2/3.$
   \end{itemize}

\begin{figure}[H]
    \centering  
    \vspace{-0.30cm} 
		\includegraphics[width=0.32\linewidth]{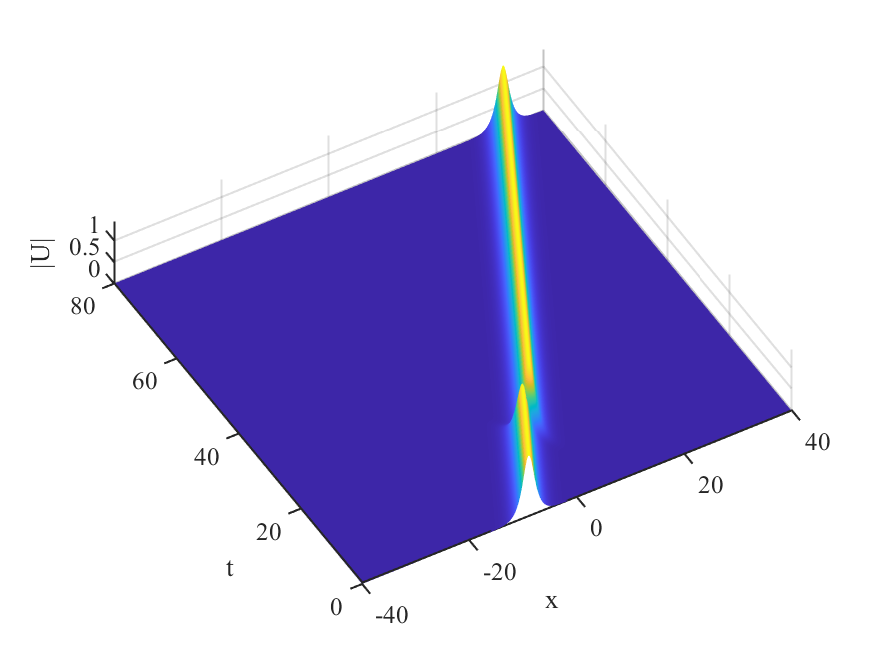}
            \includegraphics[width=0.32\linewidth]{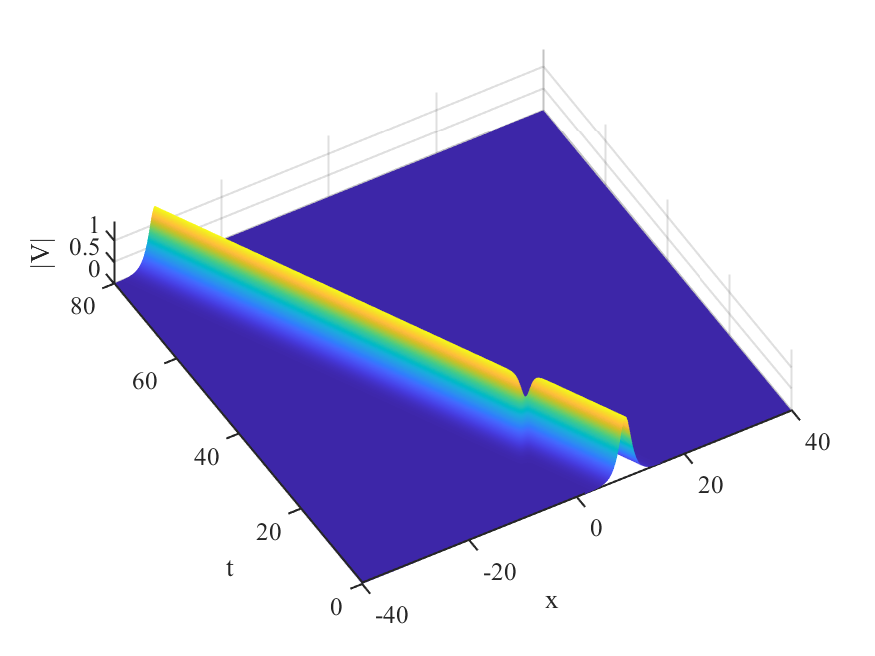}
            \includegraphics[width=0.32\linewidth]{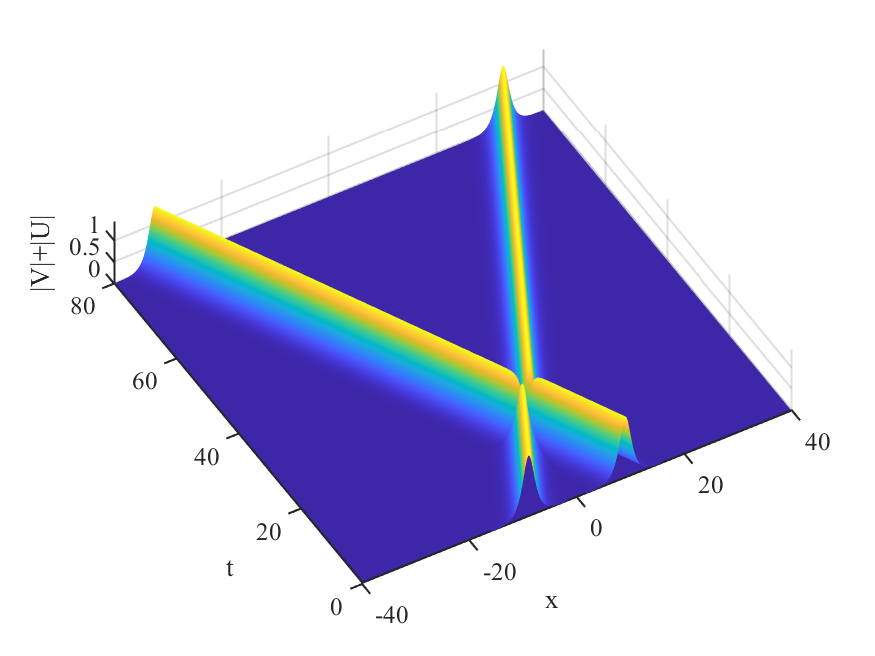}
    \caption{ Elastic  collisions of two solitons.}
    \label{fig:Elastic}
\end{figure}

\begin{figure}[H]
    \centering  
    \vspace{-0.30cm} 
		\includegraphics[width=0.32\linewidth]{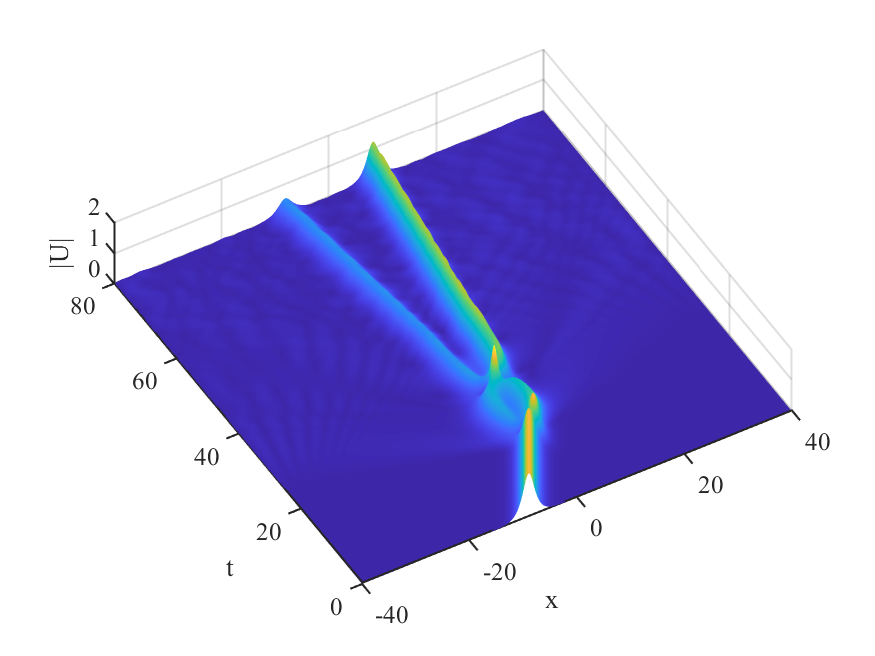}
             \includegraphics[width=0.32\linewidth]{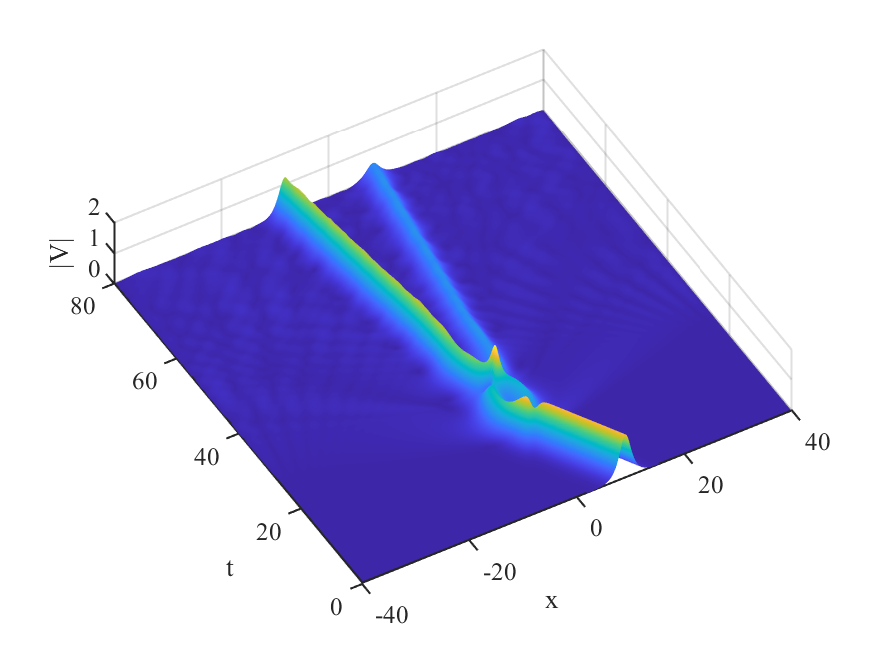}
		\includegraphics[width=0.32\linewidth]{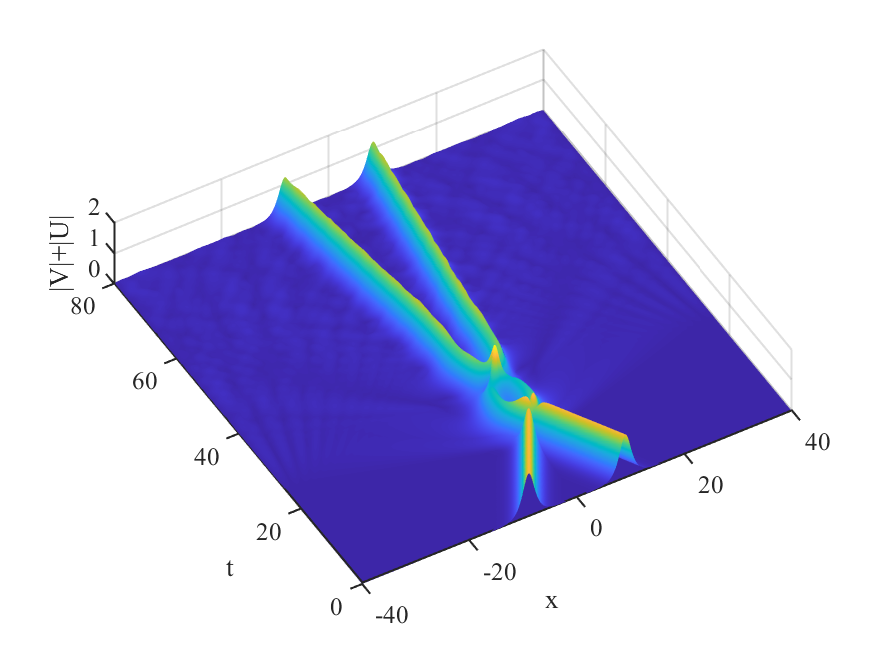}
    \caption{ Reflection of two solitons.}
    \label{fig:Reflection}
\end{figure}

\begin{figure}[H]
    \centering  
    \vspace{-0.35cm} 
	\includegraphics[width=0.32\linewidth]{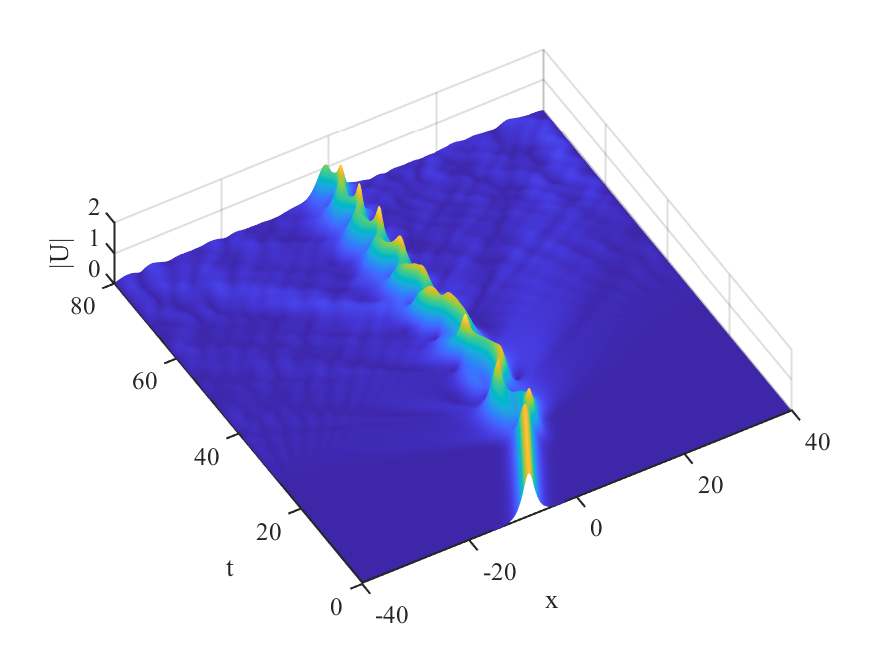}
	\includegraphics[width=0.32\linewidth]{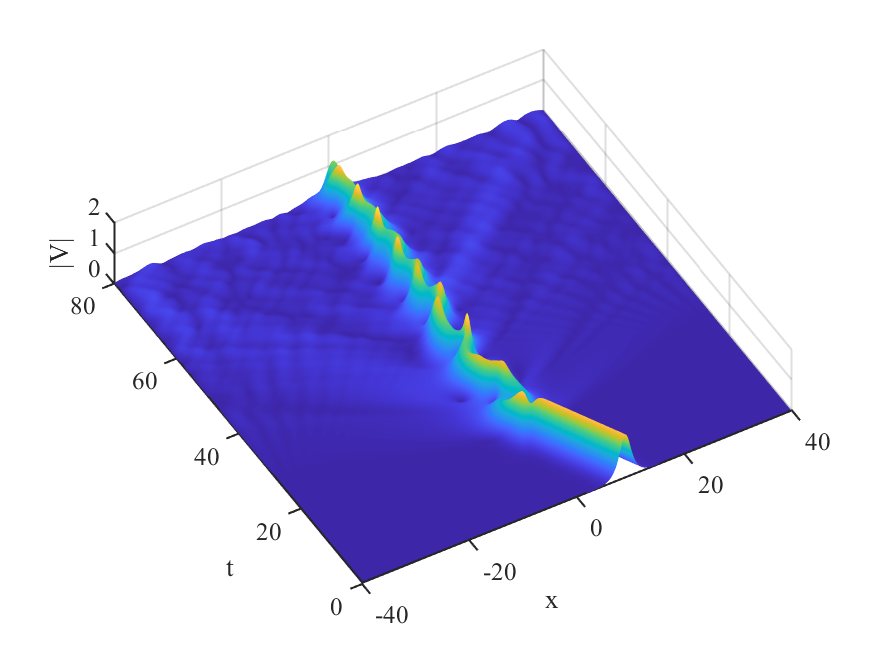}
	\includegraphics[width=0.32\linewidth]{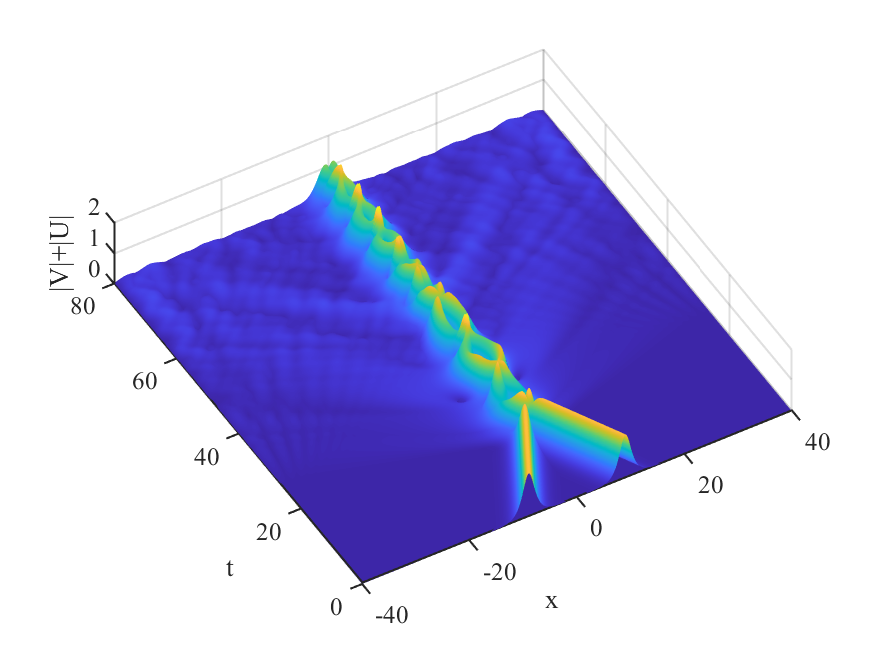}
    \caption{ Entanglement of two solitons.}
    \label{fig:Entanglement}
\end{figure}

The numerical simulation results are shown in Figures \ref{fig:Elastic}--\ref{fig:Entanglement}. We can observe the elastic collision of two solitons in Figure \ref{fig:Elastic}, in which, before collision, the solitons propagate stably with preserved amplitudes and velocities; during collision, their waveforms briefly change due to the coupling effect of the equations; and after collision, both solitons recover their original shapes and speeds, confirming perfect elasticity. In contrast, Figures \ref{fig:Reflection}--\ref{fig:Entanglement} show two different inelastic collisions, respectively. As shown in Figure \ref{fig:Reflection}, in the case of reflection, two solitons approach each other symmetrically, and after interaction, they reflect at a lower velocity with less radiation loss. Besides, it can be seen from Figure \ref{fig:Entanglement} that the solitons form oscillatory bound states after collision, exhibiting sustained mutual entanglement. All numerical simulations accurately capture these delicate nonlinear effects and align with the theoretical expectation, which verify the accuracy of the developed LRCD scheme in maintaining soliton dynamics.

\section{Conclusions}\label{sec:con}
In this paper, we developed a linear and decoupled relaxation high-order compact scheme for the CNLS system. By introducing two auxiliary relaxation variables and using a time-staggered mesh framework, the proposed scheme achieves linearization and fully decoupling, significantly improving computational efficiency compared to traditional fully implicit FD schemes (see, Example 5.2 of Section 5 in \cite{ZHOU2025}). Rigorous numerical analysis demonstrates that the scheme preserves key physical invariants (e.g., mass and energy) at the discrete level, see Theorem \ref{thm:conservation}. Most importantly, for the first time, we established the unconditional and optimal convergence analysis for the primal variables in the discrete $H^1$-norm and for the relaxation variables in the discrete $L^2$-norm, without any coupling conditions on discretization mesh parameters in multi-dimensional space situations (see, Theorem \ref{thm:Convergence}). The analysis framework developed here can also be extended to derive the unconditional optimal-order error estimates for the relaxation method applied to other NLS-type equations, for example, the Schr\"{o}dinger-Poisson equation \cite{AA2023,LHN2023}. Meanwhile, the analysis method does not have specific requirements for the spatial discretization methods. For example, in the finite element discretization framework \cite{AA2023,LHN2023}, the unconditional optimal-order error estimates can also be derived.

Finally, we notice that the convergence order in the $L^\infty$-norm is suboptimal and can only guarantee the $L^\infty$-norm boundedness of the numerical solutions. Therefore, an interesting task could be the development of an unconditional optimal-order error analysis in the $L^\infty$-norm.




\section*{Declarations}
\begin{itemize}
   
    \item \textbf{Funding}~   This work was supported in part by the National Natural Science Foundation of China (No. 12131014), by the Shandong Provincial Natural Science Foundation (No. ZR2024MA023), and by the Fundamental Research Funds for the Central Universities (No. 202264006).
   
    
       \item  \textbf{Data Availability}~  Data will be made available on request.
    
        \item \textbf{Conflict of Interest}~  The authors declare no competing interests.	

    \item \textbf{Author Contributions} {Y. Gao}: Methodology, Formal analysis, Software, Writing- Original draft. {H. Fu}: Conceptualization, Supervision, Methodology, Writing- Reviewing and Editing, Funding acquisition. {X. Wang}: Methodology, Formal analysis,  Writing- Original draft.
\end{itemize}

\bibliographystyle{spmpsci}
\bibliography{Ref_CNLS} 

\end{document}